\documentclass[12pt]{article}
\usepackage{times,amsmath,amsfonts,amsthm,amssymb,eucal}
\usepackage{algorithm}
\usepackage{algorithmic}
\usepackage{graphicx,ifthen}
\usepackage{wrapfig}
\usepackage{latexsym}
\usepackage{color}
\usepackage{mathrsfs}
\usepackage{bm}
\usepackage{bbm}
\usepackage{srctex}

\newtheorem{theorem}{Theorem}%[section]
\newtheorem{lemma}{Lemma}%[section]

\newtheorem{remark}{Remark}%[section]
\newtheorem{corollary}{Corollary}%[section]
\newtheorem{proposition}{Proposition}%[section]

\DeclareGraphicsExtensions{.gif,.pdf,.png,.jpg,.tiff}
\begin{document}
\renewcommand{\baselinestretch}{1.54}
\def\Quote{\begin{quotation}\normalfont\small}
\def\EndQuote{\end{quotation}\rm}
\def\BigHeading{\bfseries\Large}\def\MediumHeading{\bfseries\large}
\def\bct{\begin{center}}
\def\ect{\end{center}}
\font\BigCaps=cmcsc9 scaled \magstep 1
\font\BigSlant=cmsl10    scaled \magstep 1
\def\lbk{\linebreak}
\def\Report{Association Study and
Misspecified Mixed Model Analysis}
\def\Author{Jiang, Li, Paul, Yang, and Zhao}
\pagestyle{myheadings} \markboth{\Author}{\Report} \thispagestyle{empty}
\bct{\BigHeading High-dimensional Genome-wide Association Study\\
and Misspecified Mixed Model Analysis}\\\vskip10pt
\BigCaps Jiming ${\rm Jiang}^{\dagger}$,
Cong ${\rm Li}^{\ddagger}$,
Debashis ${\rm Paul}^{\dagger}$,
Can ${\rm Yang}^{\ddagger}$,
and Hongyu ${\rm Zhao}^{\ddagger}$\lbk
\BigSlant $\dagger$University of California, Davis and $\ddagger$Yale
University
\ect

\begin{abstract}
We study behavior of the restricted maximum likelihood (REML)
estimator under a misspecified linear mixed model (LMM) that has
received much attention in recent gnome-wide association studies.
The asymptotic analysis establishes consistency of the REML
estimator of the variance of the errors in the LMM, and convergence
in probability of the REML estimator of the variance of the random
effects in the LMM to a certain limit, which is equal to the true
variance of the random effects multiplied by the limiting proportion of
the nonzero random effects present in the LMM. The aymptotic results
also establish convergence rate (in probability) of the REML estimators
as well as a result regarding convergence of the asymptotic conditional
variance of the REML estimator. The asymptotic results are fully
supported by the results of empirical studies, which include extensive
simulation studies that compare the performance of the REML estimator
(under the misspecified LMM) with other existing methods.
\end{abstract}

\Quote
\vskip-5pt

\EndQuote

\section{Introduction}\label{sec:intro}
\hspace{4mm}
Genome-wide association study (GWAS), which typically refers to
examination of associations between up to millions of genetic
variants in the genome and certain traits of interest among
unrelated individuals, has been very successful for detecting
genetic variants that affect complex human traits/diseases in the past
eight years. According to the web resource of GWAS catalog (Hindorff
{\it et al.} 2009; {\bf http://www.genome.gov/gwastudies}), as of October,
2013, more than 11,000 single-nucleotide polymorphisms (SNPs) have been
reported to be associated with at least one trait/disease at the
genome-wide significance level  ($p$-value$\leq 5\times 10^{-8}$),
many of which have been validated/replicated in further studies. However,
these significantly associated SNPs only account for a small portion of
the genetic factors underlying complex human traits/diseases (Manolio
{\it et al.} 2009). For example, human height is a highly heritable
trait with an estimated heritability of around 80\%, that is, 80\% of
the height variation in the population can be attributed to genetic
factors (Visscher {\it et al.} 2008). Based on large-scale GWAS, about
180 genetic loci have been reported to be significantly associated with
human height (Allen {\it et al.} 2010). However, these loci together can
explain only about 5-10\% of variation of human height (Allen {\it et al.}
2010, Manolio {\it et al.} 2009, Visscher 2008). This ``gap'' between
the total genetic variation and the variation that can be explained by
the identified genetic loci is universal among many complex human
traits/diseases and is referred to as the ``missing heritability''
(Maher 2008, Manolio 2010, Manolio {\it et al.} 2009).

One possible explanation for the missing heritability is that many
SNPs jointly affect the phenotype, while the effect of each SNP is too
weak to be detected at the genome-wide significance level. To address
this issue, Yang {\it et al.} (2010) used a linear mixed model (LMM)-based
approach to estimate the total amount of human height variance that can be
explained by all common SNPs assayed in GWAS. They showed that 45\% of the
human height variance can be explained by those SNPs, providing compelling
evidence for this explanation: A large proportion of the heritability is
not ``missing'', but rather hidden among many weak-effect SNPs. These SNPs
may require a much larger sample size to be detected. The LMM-based
approach was also applied to analyze many other complex human
traits/diseases (e.g., metabolic syndrome traits, Vattikuti {\it et al.}
2012; and psychiatric disorders, Lee {\it et al.} 2012, Cross-Disorder
Group of Psychiatric Genomics Consortium 2013) and similar results have
been observed.

Statistically, the heritability estimation based on the GWAS data can be
cast as the problem of variance component estimation in high dimensional
regression, where the response vector is the phenotypic values and the
design matrix is the standardized genotype matrix (to be detailed below).
One needs to estimate the residual variance and the variance that can be
attributed to all of the variables in the design matrix. In a typical GWAS
data set, although there may be many weak-effect SNPs (e.g., $\sim 10^3$,
Stahl {\it et al.} 2012) that are associated with the phenotype, they
are still only a small portion of the total number SNPs (e.g., $10^5 \sim
10^6$). In other words, using a statistical term, the true underlying model
is sparse. However, the LMM-based approach used by Yang {\it et al.} assumes
that the effects of all the SNPs are non-zero. It follows that the assumed
LMM is misspecified. In spite of the huge impact of its results in the
genetics community, the misspecified LMM-based approach has not yet been
rigorously justified. In this paper, we provide theoretical justification
of the misspecified LMM in high-dimensional variance component estimation
by investigating the asymptotics of the restricted maximum likelihood
(REML; e.g., Jiang 2007) estimator as both the sample size and the
dimension of the vector of random effects tend to infinity. The results of
our theoretical study imply consistency of the REML estimators of some of
the important genetic quantities, such as the heritability, in spite of
the model misspecification. We also study convergence rate and asymptotic
variance property of the REML estimator. The theoretical results are fully
supported by the results of our empirical studies. Our study not only
provides theoretical support for the recent discoveries in human genetics
made by the LMM but also, for the first time, introduces the notion of
misspecified mixed model analysis (MMMA) and its asymptotic properties.

In addition to the significant impact of variance estimation in the genetic
community, the problem of estimating the residual variance in the
high-dimensional setting has drawn much attention recently. First, the problem
is interesting in its own right, as addressed in some recent papers (Fan {\it
et al.} 2012, Reid {\it et al.} 2013). Secondly, the significance tests for the
estimated coefficients in sparse regression (Lockhart {\it et al.} 2013,
Javanmard and Montanari 2013) require an estimator of the residual variance.
Our results open another door for the variance estimation in high-dimensional
regression. From a technical standpoint, our asymptotic analysis can be seen as
an application of the celebrated random matrix theory (e.g., Bai and
Silverstein 2010).

\subsection{Misspecified LMM and REML
estimation}\label{subsec:misspecified_LMM}
\hspace{4mm}

Consider a LMM that can be expressed as
\begin{eqnarray}\label{eq:LMM_model}
y&=&X\beta+\tilde{Z}\alpha+\epsilon,
\end{eqnarray}
where $y$ is an $n\times 1$ vector of observations; $X$ is a $n\times q$
matrix of known covariates; $\beta$ is a $q\times 1$ vector of unknown
regression coefficients (the fixed effects); $\tilde{Z}=p^{-1/2}Z$,
where $Z$ is an $n\times p$ matrix whose entries are random variables.
Furthermore, $\alpha$ is a $p\times 1$ vector of random effects that is
distributed as $N(0,\sigma_{\alpha}^{2} I_{p})$, $I_{p}$ being the
$p$-dimensional identity matrix, and $\epsilon$ is an $n\times 1$
vector of errors that is distributed as $N(0,\sigma_{\epsilon}^{2})$,
and $\alpha$, $\epsilon$, and $Z$ are independent. The estimation of
$\sigma_{\epsilon}^{2}$ is of main interest. Without loss of generality,
assume that $X$ is full rank.

The LMM (\ref{eq:LMM_model}) is what we call assumed model. In reality,
however, only a subset of the random effects are nonzero. More specifically, we
have $\alpha=\{\alpha_{(1)}^{'},0'\}'$, where $\alpha_{(1)}$ is the vector of
the first $m$ components of $\alpha$ ($1\leq m\leq p$), and $0$ is the
$(p-m)\times 1$ vector of zeros. Correspondingly, we have $\tilde
{Z}=[\tilde{Z}_{(1)}\;\tilde{Z}_{(2)}]$, where $\tilde{Z}_{(j)}=p^{-1/2}
Z_{(j)}, j=1,2$, $Z_{(1)}$ is $n\times m$, and $Z_{(2)}$ is $n\times (p-m)$.
Therefore, the true LMM can be expressed as
\begin{eqnarray}\label{eq:LMM_true_model}
y&=&X\beta+\tilde{Z}_{(1)}\alpha_{(1)}+\epsilon.
\end{eqnarray}

With respect to the true model (\ref{eq:LMM_true_model}), the assumed model
(\ref{eq:LMM_model}) is misspecified. We shall call the latter a misspecified
LMM, or mis-LMM. However, this may not be known to the investigator, who would
proceed with the standard mixed model analysis (e.g., Jiang 2007, ch. 1) to
obtain estimates of the model parameters, based on (\ref{eq:LMM_model}). This
is what we referred to as MMMA. In this paper, we will be focusing on REML
method (e.g., Jiang 2007, sec. 1.3.2). Furthermore, following Jiang (1996), we
consider estimation of $\sigma_{\epsilon}^{2}$ and the ratio $\gamma=
\sigma_{\alpha}^{2}/\sigma_{\epsilon}^{2}$. According to Jiang (2007, sec.
1.3.2), the REML estimator of $\gamma$, denoted by $\hat{\gamma}$, is the
solution to the equation
\begin{eqnarray}\label{eq:gamma_hat_equation}
\frac{y'P_{\gamma}\tilde{Z}\tilde{Z}'P_{\gamma}y}{{\rm tr}(P_{\gamma}
\tilde{Z}\tilde{Z}')}&=&\frac{y'P_{\gamma}^{2}y}{{\rm tr}(P_{\gamma})},
\end{eqnarray}
where
$P_{\gamma}=V_{\gamma}^{-1}-V_{\gamma}^{-1}X(X'V_{\gamma}^{-1}X)^{-1}X'V_{\gamma}^{-1}$
with $V_{\gamma}=I_{n}+\gamma\tilde{Z}\tilde{Z}'$. Equation
(\ref{eq:gamma_hat_equation}) is combined with another REML equation, which can
be expressed as
\begin{eqnarray}\label{eq:sigma_epsilon}
\sigma_{\epsilon}^{2}&=&\frac{y'P_{\gamma}^{2}y}{{\rm tr}(P_{\gamma})},
\end{eqnarray}
to obtain the REML estimator of $\sigma_{\epsilon}^{2}$, namely, $\hat
{\sigma}_{\epsilon}^{2}=y'P_{\hat{\gamma}}^{2}y/{\rm tr}(P_{\hat
{\gamma}})$.

In the context of mixed effects models, asymptotic behavior of the REML
estimators is well established (Das 1979, Cressie and Lahiri 1993, Richardson
and Welsh 1994, Jiang 1996). Note that the standard LMM is a conditional model,
on the $X$ and $Z$; hence, in particular, the matrix $Z$ is nonrandom. However,
this difference is relatively trivial. A more important difference is, as
noted, that the LMM (\ref{eq:LMM_model}) is misspecified. Nevertheless, what
appears to be striking is that the estimator $\hat{\sigma}_{\epsilon}^{2}$ is,
still, consistent. On the other hand, the estimator $\hat{\gamma}$ converges in
probability to a constant limit, although the limit may not be the true
$\gamma$. In spite of the inconsistency of $\hat{\gamma}$, when it comes to
estimating some important quantities of genetic interest, such as the
heritability (see below), REML still provides the right answer. Before
presenting any theoretical results, we first illustrate with a numerical
example that also highlights the practical relevance of our theoretical study.

\subsection{A numerical illustration}\label{subsec:numerical_SNP}
\hspace{4mm}

\begin{figure}\label{fig:heritability_REML}
\centering
\includegraphics[width=.68\linewidth]{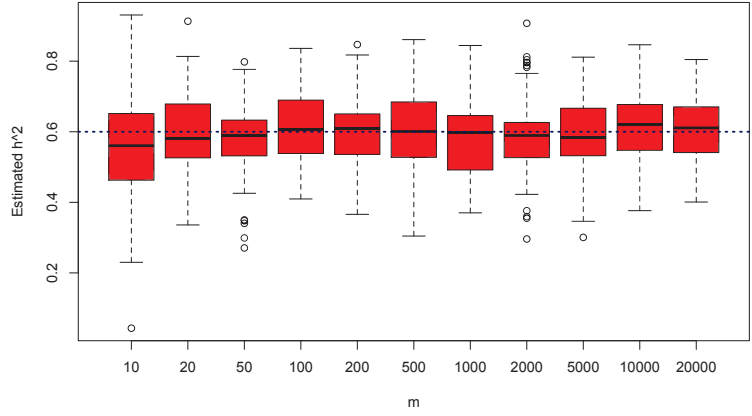}
\caption{Heritability--REML provide right answer despite
model misspecification}
\end{figure}

In GWAS, SNPs are high-density bi-allelic genetic markers. Loosely speaking,
each SNP can be considered as a binomial random variable with two trials and
the probability  of ``success'' is defined as ``allele frequency'' in genetics.
Accordingly, the genotype for each SNP can be coded as either 0, 1 or 2. In our
simulation, we first simulate the allele frequencies for $p$ SNPs, $\{f_1, f_2,
\ldots, f_p\}$, from the ${\rm Uniform}[0.05, 0.5]$ distribution, where $f_j$
is the allele frequency of the $j$-th SNP. We then simulate the genotype matrix
$U\in\{0, 1, 2\}^{n \times p}$, with rows corresponding to the
sample/individual and columns corresponding the SNP. Specifically, for the
$j$-th SNP, the genotype value of each individual is sampled from $\{0, 1, 2\}$
according to probabilities $(1-f_j)^2$, $2f_j(1-f_j)$, and $f^2_j$,
respectively. After that, each column of $U$ is standardized to have zero mean
and unit variance, and the standardized genotype matrix is denoted as $Z$. Let
$\tilde{Z}=p^{-1/2}Z$. In Yang {\it et al.} (2010), an LMM was used to describe
the relationship between a phenotypic vector $y$ and the standardized genotype
matrix $\tilde{Z}$:
\begin{eqnarray}\label{eq:phenotype_model}
y=1_{n}\mu+\tilde{Z}\alpha+\epsilon,\;\;\alpha\sim N(0,
\sigma^2_{\alpha}I_p), \epsilon\sim N(0, \sigma^2_{\epsilon}I_n),
\end{eqnarray}
where $1_{n}$ is the $n\times 1$ vector of 1's, $\mu$ is an intercept, $\alpha$
is the vector of random effects, $I_n$ is the $n\times n$ identity matrix, and
$\epsilon$ is the vector of errors. An important quantity in genetics is
``heritability'', defined as the proportion of phenotypic variance explained by
all genetic factors. For convenience, we assume that all of the genetic factors
have been captured by the SNPs in GWAS. Under this assumption, the heritability
can be characterized via the variance components in model
(\ref{eq:phenotype_model}):
\begin{eqnarray}\label{eq:h_square}
h^2&=&\frac{\sigma^2_{\alpha}}{\sigma^2_{\alpha}+\sigma^2_{\epsilon}}.
\end{eqnarray}
Note that the definition of heritability by (\ref{eq:h_square}) assumes that
$\alpha_j \sim N(0, \sigma^2_{\alpha})$ for all $j \in \{0, 1, 2, \ldots, p\}$.
However, in reality, only a subset of the SNPs are associated with the
phenotype. A correct model therefore is
\begin{eqnarray}\label{eq:phenotype_true_model}
y=1_{n}\mu+\tilde{Z}_{(1)}\alpha_{(1)}+\epsilon,\;\;\alpha \sim N(0,
\sigma^2_{\alpha}I_m),\;\;\epsilon\sim N(0, \sigma^2_{\epsilon}I_n),
\end{eqnarray}
where $m$ is the total number of SNPs that are associated with the phenotype,
$\alpha_{(1)}$ is the subvector of $\alpha$ corresponding to the nonzero
components that are associated with the SNPs, and $\tilde{Z}_{(1)}=
p^{-1/2}Z_{(1)}$, $Z_{(1)}$ being the submatrix of $Z$ corresponding to the
associated SNPs. In this case, the heritability should instead be given by
\begin{eqnarray}\label{eq:h_true}
h^2_{\rm true}&=&\frac{(m/p)\sigma^2_{\alpha}}{(m/p)\sigma^2_{\alpha}
+\sigma^2_{\epsilon}}.
\end{eqnarray}
In practice, it is impossible to identify all of the $m$ SNPs due to the
limited sample size. Therefore, we follow model (\ref{eq:phenotype_true_model})
while simulating the phenotypic values, but pretend that we do not know which
SNPs are associated with the phenotype. This means that we simply use all the
SNPs in $Z$ to estimate the variance components, ${\sigma}^2_{\alpha}$ and
${\sigma}^2_{\epsilon}$ in model (\ref{eq:phenotype_model}). The estimated
heritability is then obtained as
\begin{eqnarray}\label{eq:h_hat_square}
\hat{h}^2&=&\frac{\hat{\sigma}^2_{\alpha}}{\hat{\sigma}^2_{\alpha}
+\hat{\sigma}^2_{\epsilon}}.
\end{eqnarray}
In this illustrative simulation, we fixed $n = 2,000$, $p = 20,000$,
$\sigma^2_{\epsilon}=0.4$ and varied $m$ from $10$ to $20,000$. We also set the
variance component $\sigma^2_{\alpha}=0.6p/m$ so that the proportion of
phenotypic variance explained by genetic factors $h^2_{\rm true}=0.6$, based on
(\ref{eq:h_true}). We repeated the simulation 100 times. As shown in Figure
\ref{fig:heritability_REML}, there is almost no bias in the estimated $h^2$
regardless of the underlying true model, whether it is sparse (i.e., $m/p$ is
close to zero) or dense (i.e., $m/p$ is close to one). This suggests that the
REML works well in providing unbiased estimator of the heritability despite the
model misspecification.

\subsection{Outline of theoretical results}\label{subsec:outline_theoretical}
\hspace{4mm}

Throughout this paper, we assume that $q$, the dimension of $\beta$,
is fixed, while $n$, $p$, and $m$ increase. For the simplicity of
illustration, let us first assume that $n,p,m\rightarrow\infty$ such that
\begin{eqnarray}\label{eq:m_n_p_condition}
\frac{n}{p}\longrightarrow\tau,\;\;\;\frac{m}{p}\longrightarrow\omega,
\end{eqnarray}
where $0<\tau,\omega\leq 1$ are constants. Note that $\tau$ is the
limiting ratio of the sample size and the number of random effects,
while $\omega$ is the limiting proportion of the nonzero random effects.
First consider the case where the entries of $Z$ are i.i.d. The point is
that the more realistic case where the entries of $Z$ are standardized
(see below) can be handled by utilizing the results for the i.i.d. case,
and some inequalities on the difference, or perturbation (see below),
between the two cases.

Suppose that the true variance components, $\sigma_{\alpha}^{2},
\sigma_{\epsilon}^{2}$ are positive, and (\ref{eq:m_n_p_condition}) holds.
Then, (i) with probability tending to one, there is a REML estimator,
$\hat{\gamma}$, such that $\hat{\gamma}\stackrel{\rm
P}{\longrightarrow}\omega\gamma_{0}$, where $\gamma_{0}$ is the true $\gamma$;
(ii) $\hat{\sigma}_{\epsilon}^{2} \stackrel{\rm
P}{\longrightarrow}\sigma_{\epsilon 0}^{2}$, where $\hat
{\sigma}_{\epsilon}^{2}$ is the REML estimator given by
(\ref{eq:sigma_epsilon}) with $\gamma=\hat{\gamma}$, as in (i), and
$\sigma_{\epsilon 0}^{2}$ is the true $\sigma_{\epsilon}^{2}$.

As far as the consistency is concerned, condition (\ref{eq:m_n_p_condition})
can be relaxed to
\begin{eqnarray}\label{eq:m_n_p_condition_relaxed}
\liminf\left(\frac{m\wedge n}{p}\right)>0,\;\;
\limsup\left(\frac{m\vee n}{p}\right)\leq 1.
\end{eqnarray}
so that, with probability tending to one, that there exist REML
estimators, $\hat{\gamma}, \hat{\sigma}_{\epsilon}^{2}$, such that
(i) $\hat{\sigma}_{\epsilon}^{2}\stackrel{\rm
P}{\longrightarrow}\sigma_{\epsilon 0}^{2}$, in other words, the REML
estimator of $\sigma_{\epsilon}^{2}$ is consistent; and
(ii) the adjusted REML estimator of $\gamma$ is consistent, that is,
$(p/m)\hat{\gamma}\stackrel{\rm P}{\longrightarrow}\gamma_{0}$.

{\bf Note.} The latest asymptotic result may explain what has been observed in
Figure 1. %\ref{fig:heritability_REML}
Note that the estimated heritability, (\ref{eq:h_hat_square}), can be written
as
\begin{eqnarray}\label{eq:h_hat_square_expr}
\hat{h}^{2}&=&\frac{(m/p)(p/m)\hat{\gamma}}{1+(m/p)(p/m)\hat{\gamma}}.
\end{eqnarray}
On the other hand, the true heritability, (\ref{eq:h_true}), can be written as
\begin{eqnarray}\label{eq:h_true_expr}
h_{\rm true}^{2}&=&\frac{(m/p)\gamma_{0}}{1+(m/p)\gamma_{0}}.
\end{eqnarray}
Because $(p/m)\hat{\gamma}$ converges in probability to $\gamma_{0}$, when we
replace the $(p/m)\hat{\gamma}$ in (\ref{eq:h_hat_square_expr}) by
$\gamma_{0}$, the resulting first-order approximation of
(\ref{eq:h_hat_square_expr}) is exactly (\ref{eq:h_true_expr}). It should also
be noted that condition (\ref{eq:m_n_p_condition_relaxed}) requires that the
limiting lower bound be positive. This may explain why the bias for $m=10$ in
Figure 1 %\ref{fig:heritability_REML}
is much more significant compared to other
cases, because the ratio $m/p$ in this case, $10/20000=0.0005$, is fairly close
to zero.

As mentioned, the asymptotic results can be extended to the case where the
design matrix, $Z$, for the random effects is standardized. Let
$U=(u_{ik})_{1\leq i\leq n, 1\leq k\leq p}$ whose entries are i.i.d. Define
$Z=(z_{ik})_{1\leq i\leq n,1\leq k\leq p}$, where $z_{ik}
=(u_{ik}-\bar{u}_{k})/s_{k}$ with $u_{k}=n^{-1}\sum_{i=1}^{n}u_{ik}$ and
$s_{k}^{2}=(n-1)\sum_{i=1}^{n}(u_{ik}-\bar{u}_{k})^{2}$. In other words, the
new $Z$ matrix has the sample mean equal to $0$ and sample variance equal to
$1$ for each column. We then define $\tilde{Z}=p^{-1}Z$, and proceed as in
(\ref{eq:LMM_model}). Also, as noted, in GWAS, the entries of $U$ are generated
from a discrete distribution which assigns the probabilities $\theta^{2},
2\theta(1-\theta), (1-\theta)^{2}$ to the values $0, 1, 2$, where $\theta$ is
pre-specified so that $\theta\in (0.05,0.5)$; however, there is also interest
in the case where the entries of $U$ are normal. Under the discrete
distribution, it makes no difference if we standardize the discrete
distribution so that is has mean $0$ and variance $1$, so, without loss of
generality, the entries of $U$ are $u_{ik}=(d_{ik} -\mu)/\sigma$, where
$d_{ij}$ has the above discrete distribution, $\mu={\rm
E}(d_{ik})=2(1-\theta)$, and $\sigma^{2}={\rm var}(u_{ik}) =2\theta(1-\theta)$.

Both the Gaussian and discrete cases can be treated under the
framework of the following broader class of distributions (e.g., Hsu
{\it et al.} 2012). Let $\xi_{1},\dots,\xi_{n}$ be random variables. We
say $\xi=(\xi_{i})_{1\leq i\leq n}$ is sub-Gaussian if there exists
$\sigma>0$ such that for all $\lambda\in R^{n}$ we have ${\rm
E}(e^{\lambda'\xi})\leq e^{|\lambda|^{2}\sigma^{2}/2}$. The asymptotic
results regarding the MMMA are extended to the sub-Gaussian class.

In addition to the consistency results, we also study convergence
rate and asymptotic variance property of the REML estimator under the
mis-LMM. The results provide further insights into the asymptotic
behavior of these estimators.

\section{Preliminaries}\label{sec:preliminaries}
\hspace{4mm}

A key component for our proofs is the following celebrated result in random
matrix theory (e.g., Paul and Aue 2013). Let $Z$ be an $n\times p$ matrix whose
entries are i.i.d., complex-valued random variables with mean $0$ and variance
$1$, where $n\rightarrow\infty$ as $p\rightarrow\infty$ such that
$n/p\rightarrow\tau$, as in (\ref{eq:m_n_p_condition}). We are interested in
the asymptotic behavior of the empirical spectral distribution (ESD) of
$S=p^{-1}ZZ'$, defined as
\begin{eqnarray*}
F^{S}(x)&=&\frac{1}{n}\sum_{k=1}^{n}1_{(\lambda_{k}\leq x)},\;\;x\in
\mathbb{R},
\end{eqnarray*}
where $\lambda_{1},\dots,\lambda_{n}$ are the eigenvalues of $S$.

\begin{lemma}\label{lemma:MP_law}
(Mar\v{c}enko-Pastur law) Suppose (\ref{eq:m_n_p_condition}) holds. Then, as
$p\rightarrow\infty$, the ESD of $S$ converges almost surely (a.s.) in
distribution to the Mar\v{c}enko-Pastur (M-P) law, $F_{\tau}$, whose p.d.f. is
given by
\begin{eqnarray*}
f_{\tau}(x)&=&\frac{1}{2\pi\tau x}\sqrt{\{b_{+}(\tau)-x\}\{x-b_{-}(\tau)
\}},
\end{eqnarray*}
if $b_{-}(\tau)\leq x\leq b_{+}(\tau)$, and $f_{\tau}(x)=0$ elsewhere,
where $b_{\pm}(\tau)=(1\pm\sqrt{\tau})^{2}$.
\end{lemma}

A result that is frequently referred to is the following corollary of Lemma
\ref{lemma:MP_law}, which is a consequence of convergence in distribution
(e.g., Jiang 2010, p. 45).

\begin{corollary}\label{cor:MP_corollary}
Under the assumptions of Lemma \ref{lemma:MP_law}, we have, for any positive
integer $l$, $n^{-1}{\rm tr}(S^{l})\stackrel{\rm
a.s.}{\longrightarrow}\int_{b_{-}(\tau)}^{b_{+}(\tau)}x^{l}f_{\tau}(x)dx$ as
$p\rightarrow\infty$.
\end{corollary}

The next result is regarding the extreme eigenvalues of $S$ (e.g., Bai 1999,
th. 2.16). Let $\lambda_{\min}(S)$ (respectively, $\lambda_{\max}(S)$) denote
the smallest (largest) eigenvalues of $S$.

\begin{lemma}\label{lemma:extreme_eigenvalues}
Suppose that, in addition to the assumptions of Lemma \ref{lemma:MP_law}, the
fourth moment of the entries of $Z$ are finite. Then, we have, as
$p\rightarrow\infty$, $\lambda_{\min}(S)\stackrel{\rm a.s.}{\longrightarrow}
b_{-}(\tau)$ and $\lambda_{\max}(S)\stackrel{\rm a.s.}{\longrightarrow}
b_{+}(\tau)$.
\end{lemma}

Let $\xi_{1},\dots,\xi_{n}$ be random variables. We say $\xi=(\xi_{i})_{1 \leq
i\leq n}$ is sub-Gaussian if there exists $\sigma>0$ such that for all
$\lambda\in \mathbb{R}^{n}$ we have ${\rm E}(e^{\lambda'\xi})\leq
e^{|\lambda|^{2}\sigma^{2}/2}$. The Gaussian distribution, of course, is a
member of the sub-Gaussian class. The following is a restatement of Lemma 5.5
of Vershynin (2011).

\begin{lemma}\label{lemma:subGaussian}
A random variable $\xi$ is sub-Gaussian if any of the following equivalent
conditions hold:
\begin{itemize}
\item[(I)] ${\rm E}(e^{\xi^2/K_1^2})<\infty$ for some $0<K_1<\infty$;

\item[(II)] $\{{\rm E}(|\xi|^q)\}^{1/q}\leq K_2\sqrt{q}$ for all $q\geq 1$, for
some $0<K_2<\infty$.

If, moreover, ${\rm E}(\xi)=0$, then the following is equivalent to
(I) and (II):

\item[(III)] ${\rm E}(e^{t\xi})\leq e^{t^2 K_3^2}$ for all $t\in R$, for some $0
<K_3<\infty$.
\end{itemize}
\end{lemma}

Define the \textit{sub-Gaussian norm} of a random variable $\xi$ as
\begin{eqnarray*}
\|\xi\|_{\psi_2}&\equiv&\sup_{q\geq 1}\left\{q^{-1/2}({\rm E}|\xi|^q)^{1/q}
\right\}.
\end{eqnarray*}
Clearly, by (II) of Lemma \ref{lemma:subGaussian}, $\xi$ is a sub-Gaussian
random variable if and only if $|\xi|_{\psi_2}<\infty$. One of the useful
characteristics of sub-Gaussianity is that it is preserved under linear
combinations. Specifically, we have the following result.

\begin{lemma}\label{lemma:subGaussian_linear_norm}
(Vershynin 2011, lem. 5.9). Suppose that $X_1,\ldots, X_n$ are independent
sub-Gaussian random variables, and $b_1,\ldots,b_n \in \mathbb{R}$ are
nonrandom. Then $\sum_{i=1}^n b_i X_i$ is sub-Gaussian and, for some $C > 0$,
we have
\begin{eqnarray*}
\left\|\sum_{i=1}^n b_i X_i \right\|_{\psi_2}^2&\leq&C\sum_{i=1}^n
b_i^2\|X_i\|_{\psi_2}^2.
\end{eqnarray*}
\end{lemma}

Lemma \ref{lemma:subGaussian_linear_norm}  follows easily from the equivalent
characterizations in Lemma \ref{lemma:subGaussian}, specifically, by using the
moment generating function. The following simple corollary is very useful for
our applications.

\begin{corollary}\label{cor:subGaussian_linear_norm}
Let $X_1,\ldots,X_n$ be independent with $\max_{1\leq i \leq
n}\|X_i\|_{\psi_2}\leq K<\infty$. Then $\sum_{i=1}^n b_i X_i$ is sub-Gaussian
and, for some $C > 0$, we have
\begin{eqnarray*}
\left\|\sum_{i=1}^n b_i X_i\right\|_{\psi_2}^2&\leq&C K^2
\left(\sum_{i=1}^n b_i^2\right).
\end{eqnarray*}
\end{corollary}

The following result, due to Rudelson and Vershynin (2013), is a
concentration inequality for quadratic forms involving a random
vector with independent sub-Gaussian components. It is referred to
as \textit{Hanson-Wright inequality}. For any matrix $A$ of real
entries, the spectral norm of $A$ is defined as $\|A\|=
\lambda_{\max}^{1/2}(A'A)$ and the Euclidean norm is defined as
$\|A\|_{2}={\rm tr}^{1/2}(A'A)$.

\begin{proposition}\label{prop:Hanson_Wright}
Let $\boldsymbol{\xi} = (\xi_1,\ldots,\xi_n)'$, where the $\xi_i$'s are
independent random variables satisfying ${\rm E}(\xi_i) = 0$ and $\max_{1\leq i
\leq n}\|\xi_i\|_{\psi_2}\leq K<\infty$. Let $A$ be an $n\times n$ matrix.
Then, for some constant $c>0$, we have, for any $t>0$,
\begin{eqnarray*}
{\rm P}\{|\boldsymbol{\xi}'A\boldsymbol{\xi}-{\rm E}(\boldsymbol{\xi}'
A \boldsymbol{\xi})|>t\}&\leq&2\exp\left\{-c\min\left(\frac
{t^2}{K^4\|A\|_{2}^{2}},\;\frac{t}{K^2\|A\|}\right)\right\}.
\end{eqnarray*}
\end{proposition}

In the settings that we are interested in, we have ${\rm E}(\xi_i^2)=1$
for all $i$ and so ${\rm E}(\boldsymbol{\xi}'A\boldsymbol{\xi})$ reduces
to ${\rm tr}(A)$.

The next result, well known in random matrix theory (e.g., Bai and
Silverstein 2010; sec. A.5, A.6), is regarding perturbation of the ESD.

\begin{lemma}\label{lemma:ESD_comparison}
For any $n\times p$ matrices $A, B$ we have
\begin{itemize}
\item[(i)] $\|F^{AA'}-F^{BB'}\|\leq n^{-1}{\rm rank}(A-B)$, where for a
real-valued function $g$ on $\mathbb{R}$, $\|g\|=\sup_{x\in \mathbb{R}}|g(x)|$;

\item[(ii)] $L^{4}(F^{AA'},F^{BB'})\leq 2n^{-2}(\|A\|_{2}^{2}+\|B\|_{2}^{2})
\|A-B\|_{2}^{2}$, where the Levy distance between two distributions, $F$ and
$G$ on $\mathbb{R}$, is defined as $L(F,G)=\inf\{\epsilon>0: F(x-\epsilon)
-\epsilon\leq G(x)\leq F(x+\epsilon)+\epsilon\}$.
\end{itemize}
\end{lemma}

The following result is implied by Lemma 2 of Bai and Yin (1993).

\begin{lemma}\label{lemma:columnwise_mean_dev}
Suppose that $X_{ij}, i,j=1,2,\dots$ are i.i.d. with ${\rm
E}(X_{11}^{2})<\infty$. Then, we have $\max_{1\leq j\leq
n}\left|\bar{X}_{j}-{\rm E}(X_{11})\right|\stackrel {\rm
a.s.}{\longrightarrow}0$, where $\bar{X}_{j}=n^{-1}\sum_{i=1}^{n}X_{ij}$.
\end{lemma}

Lemma \ref{lemma:ESD_comparison} and Lemma \ref{lemma:columnwise_mean_dev} are
used to study the asymptotic ESD of symmetric random matrices involving the
standardized design matrix. Note that the standardized design matrix can be
expressed as $Z=(U-\bar{u}\otimes 1_{n})D_{s}^{-1}$, where
$\bar{u}=(\bar{u}_{1},\dots,\bar{u}_{p})$, and $D_{s}={\rm
diag}(s_{1},\dots,s_{p})$ (where $\otimes$ denotes the Kronecker product). Let
$A$ be the matrix associated with the REML estimation (see the beginning of the
proof of Theorem \ref{thm:consistency_REML} below). Consider $\Psi=p^{-1}
\zeta\zeta'$, where $\zeta=A'Z$ and $A$ is $n\times(n-q)$ satisfying $A'X=0$
and $A'A=I_{n-q}$. The following corollary is proved in Section
\ref{sec:proofs}.

\begin{corollary}\label{cor:ESD_psi_convergence}
Under the assumptions of Lemma \ref{lemma:MP_law}, the ESD of $\Psi$ converges
a.s. in distribution to the M-P law. Furthermore, under the assumptions of
Lemma \ref{lemma:extreme_eigenvalues}, $\lambda_{\min}(\Psi)$ and
$\lambda_{\max}(\Psi)$ converge a.s. $b_{-}(\tau)$ and $b_{+}(\tau)$,
respectively.
\end{corollary}

\section{Main theoretical results}\label{sec:theory}
\hspace{4mm}

First we state a result regarding the consistency of the misspecified
REML estimator of $\sigma_{\epsilon}^{2}$, $\hat{\sigma}_{\epsilon}^{2}$,
and convergence in probability of the misspecified REML estimator of
$\gamma$, $\hat{\gamma}$. Throughout this section, the design matrix,
$Z$, is assumed to be the standardized, as described near the end of
Section 1, where the entries of $U$ are i.i.d. sub-Gaussian.

\begin{theorem}\label{thm:consistency_REML}
Suppose that the true $\sigma_{\alpha}^{2}, \sigma_{\epsilon}^{2}$ are
positive, and (\ref{eq:m_n_p_condition}) holds. Then,

\begin{itemize}
\item[(i)] With probability tending to one, there is a REML estimator,
$\hat{\gamma}$, such that $\hat{\gamma}\stackrel{\rm P}{\longrightarrow}
\omega\gamma_{0}$, where $\gamma_{0}$ is the true $\gamma$.

\item[(ii)] $\hat{\sigma}_{\epsilon}^{2}\stackrel{\rm P}{\longrightarrow}
\sigma_{\epsilon 0}^{2}$, where $\hat{\sigma}_{\epsilon}^{2}$ is (4) with
$\gamma=\hat{\gamma}$, as in (i), and $\sigma_{\epsilon 0}^{2}$ is the true
$\sigma_{\epsilon}^{2}$.
\end{itemize}
\end{theorem}

\begin{remark}\label{rem:limit_hat_gamma}
It is interesting to note that the limit of $\hat{\gamma}$ in (i) depends on
$\omega$, but not $\tau$. More specifically, the limit is equal to the true
$\gamma$ multiplied by $\omega$, the limiting proportion of the nonzero random
effects (see the remark below (\ref{eq:m_n_p_condition})). The result seems
totally intuitive.
\end{remark}

\begin{remark}
On the other hand, part (ii) of Theorem \ref{thm:consistency_REML} states that
the REML estimator of $\sigma_{\epsilon}^{2}$ is consistent in spite of the
model misspecification.
\end{remark}

As far as the consistency of $\hat{\sigma}_{\epsilon}^{2}$ is concerned,
condition (\ref{eq:m_n_p_condition}) can be relaxed. We state this as a
corollary of Theorem \ref{thm:consistency_REML}.

\begin{corollary}\label{cor:consistency_REML}
Suppose that, in Theorem \ref{thm:consistency_REML}, condition
(\ref{eq:m_n_p_condition}) is weakened to (\ref{eq:m_n_p_condition_relaxed}).
Then, with probability tending to one, there are REML estimators,
$\hat{\gamma}, \hat{\sigma}_{\epsilon}^{2}$, such that
\begin{itemize}
\item[(i)] $\hat{\sigma}_{\epsilon}^{2}\stackrel{\rm
P}{\longrightarrow}\sigma_{\epsilon 0}^{2}$, in other words, the REML estimator
of $\sigma_{\epsilon}^{2}$ is consistent;
\item[(ii)] The adjusted REML estimator of $\gamma$ is consistent, that is,
$(p/m)\hat{\gamma}\stackrel{\rm P}{\longrightarrow}\gamma_{0}$.
\end{itemize}
\end{corollary}

Another consequence of Theorem \ref{thm:consistency_REML} may be regarded as an
extension of the well-known result on consistency of the REML estimator (e.g.,
Jiang 1996), which is based on conditioning on $Z$.

\begin{corollary}\label{cor:consistency_REML_correct}
Suppose that $m=p$, that is, the LMM is correctly specified. Then, as
$n,p\rightarrow\infty$ such that (\ref{eq:m_n_p_condition_relaxed}) holds with
$m=n$, there are REML estimators $\hat{\gamma}$ and
$\hat{\sigma}_{\epsilon}^{2}$ such that $\hat{\gamma}\stackrel{\rm
P}{\longrightarrow}\gamma_{0}$ and $\hat{\sigma}_{\epsilon}^{2}\stackrel {\rm
P}{\longrightarrow}\sigma_{\epsilon 0}^{2}$; in other words, the REML
estimators are consistent without conditioning on $Z$.
\end{corollary}

Given the consistency of $\hat{\sigma}_{\epsilon}$, more precise
asymptotic behavior of the latter is of interest. As noted, the estimation
of $\sigma_{\epsilon}^{2}$ is also of main practical interest. The
following result establishes convergence rate of the REML estimator of
$\sigma_{\epsilon}^{2}$ as well as that of the adjusted REML estimator
of $\gamma$.

\begin{theorem}\label{thm:rate_REML}
If, in the assumption of Theorem \ref{thm:consistency_REML},
(\ref{eq:m_n_p_condition}) is strengthened to
\begin{eqnarray}\label{eq:m_n_p_rate}
\sqrt{n}\left|\frac{n}{p}-\tau\right|\rightarrow 0,\;\;\;\sqrt{n}\left|
\frac{m}{p}-\omega\right|\rightarrow 0,
\end{eqnarray}
then we have $\hat{\gamma}-\omega\gamma_{0}=O_{\rm P}(\sqrt{\log n/n})$
and $\hat{\sigma}_{\epsilon}^{2}-\sigma_{\epsilon 0}^{2}=O_{\rm P}(\sqrt
{\log n/n})$. More specifically, we have $\hat{\sigma}_{\epsilon}^{2}
-\sigma_{\epsilon 0}^{2}=t_{1}+t_{2}$, where $t_{1}=O_{\rm P}(\sqrt
{\log n/n})$ and $t_{2}=o_{\rm P}(\sqrt{\log n/n})$. The leading term,
$t_{1}$, has the property that its conditional variance on $Z$,
multiplied by $n$, converges in probability to a constant limit. It is
in the latter sense that the REML estimator of $\sigma_{\epsilon}^{2}$
has a convergent asymptotic conditional variance at the rate $1/n$.
\end{theorem}

The proofs of the theorems are given in Section \ref{sec:proofs}.

{\bf Note.} Although, throughout this paper, we have assumed that the dimension
of $\beta$, $q$, is fixed (see the beginning of Section
\ref{subsec:outline_theoretical}), the proofs show that the results of Theorem
\ref{thm:consistency_REML} and Theorem \ref{thm:rate_REML} remain valid as long
as $q=o(\sqrt{n})$. Another consequence of the latter condition is following.
Throughout this paper, the matrix of covariates, $X$ in (\ref{eq:LMM_model}),
is considered fixed. This is equivalent to the assumption that $X$ and $Z,
\epsilon$ are independent. However, as long as $q=o(\sqrt{n})$, the
independence of $X$ and $Z$ is asymptotically ignorable in that the results of
Theorem \ref{thm:consistency_REML} and Theorem \ref{thm:rate_REML} continue to
hold even if $X$ is not independent with $Z$. This is because the REML
procedure depends on $X$ only through the matrix $A$, which has the property
that $A'X=0$ and $A'A=I_{n-q}$. Furthermore, as argued near the end of the
proof of Theorem \ref{thm:consistency_REML} (see Section
\ref{subsec:proof_thm_consistency}), what is actually at play is the matrix
$AA'=I_{n}-P_{X}$, and $P_{X}$ has rank $q=o(\sqrt{n})$. It turns out that,
under the latter condition, $P_{X}$ is ignorable in all of our asymptotic
arguments; in other words, one can replace $AA'$ by $I_{n}$ and the results do
not change.

\section{More simulation studies}\label{sec:simulation}
\hspace{4mm}

To demonstrate our theoretical results numerically, we carry out more
comprehensive simulation study following the same procedures as described in
Section \ref{subsec:numerical_SNP}. The $h^2$ was also set at 0.6 ($\sigma_e^2
= 0.4$ and $\gamma = 1.5$). We fix the ratio $\tau=n/p=0.1$ and varied
$\omega=m/p$ from 0.001 to 1. We examine the performance of the REML, under the
mis-LMM, in estimating $\gamma$ and $\sigma^2_e$ as $n$ varies from 1000 to
5000. The performance of the adjusted REML estimator of $\gamma$ for $\omega =
0.01$ is shown in Figure 2. %\ref{fig:adjusted_REML}
It appears that the adjusted REML always gives nearly unbiased estimate of
$\gamma$, confirming our observations confirming our observations in Section
\ref{subsec:numerical_SNP} and theoretical results, namely, part (ii) of
Corollary \ref{cor:consistency_REML}. More importantly, as both $n$ and $p$
increase (with $n/p$ fixed at 0.1), the standard deviation of the estimate
decreases.
\begin{figure}\label{fig:adjusted_REML}
\centering
\includegraphics[width=.68\linewidth]{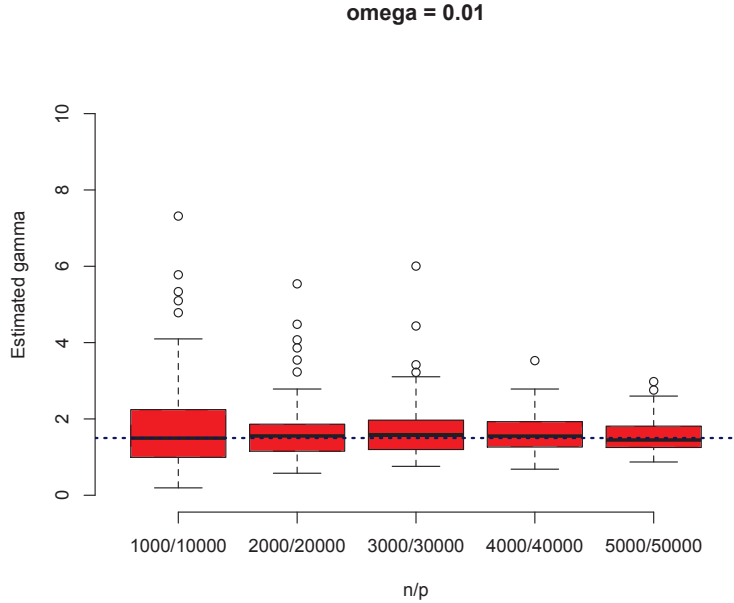}
\caption{Adjusted REML estimation of $\gamma$ for different $n$ and $p$
($\omega = 0.01$).}
\end{figure}

As noted, several other methods for high dimensional variance estimation have
been proposed recently. As a comparison, we examine the performances of two of
these methods, refitted cross validation (c.v.) (Fan {\it et al.} 2012) and
scaled lasso (Sun and Zhang 2012), in estimating $\sigma^2_e$ under the
misspecified LMM. The results for $n=2000$, $p=20000$ are shown in Figure 3.
%\ref{fig:comparsion_variance_estimators}
Again, the REML estimator appears to be unbiased regardless of the value of
$m$. On the other hand, the competing methods tend to have much larger bias,
especially when $m$ is large. This is not surprising because the competing
methods are largely based on the sparsity assumption that $m$ is relatively
small compared to $p$. Indeed, when $m = 20$, the biases and standard
deviations of the competing methods are quite small. In the latter case, the
competing method may outperform the REML in terms of mean squared error (MSE).
However, the REML performs well consistently across a much broader range of
$m$, as demonstrated by Figure 3. %\ref{fig:comparsion_variance_estimators}

\begin{figure}\label{fig:comparsion_variance_estimators}
\centering
\includegraphics[width=.68\linewidth]{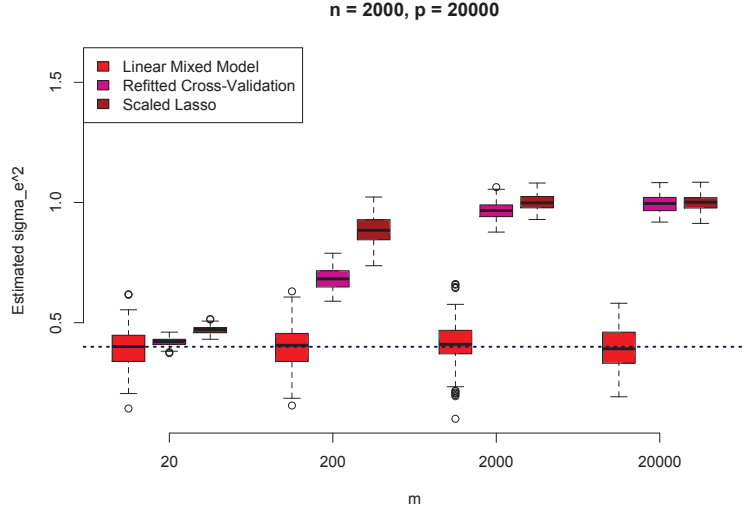}
\caption{Comparison of estimators of $\sigma_\epsilon^2$ with refitted C.V. and
scaled lasso for different $m$ ($\omega = 0.01$).}
\end{figure}

\section{Proofs}\label{sec:proofs}

\subsection{Proof of Corollary \ref{cor:ESD_psi_convergence}}\label{subsec:proof_cor_psi}
\hspace{4mm}

Note that $\zeta =(M-L)D_{s}^{-1}$, where $M=A'U$ and $L= A'\bar{u}\otimes
1_{n}$, and that $M$ is $(n-q)\times p$ whose entries are independent
sub-Gaussian, with mean $0$, variance $1$, and $A'A=I_{n -q}$. Furthermore,
write $\tilde{M}=M/\sqrt{p}$ and $\tilde{L}=L/\sqrt{p}$. By Lemma
\ref{lemma:MP_law}, the ESD of $\tilde{M}\tilde{M}'$ converges a.s. in
distribution to the M-P law. On the other hand, write
$\tilde{B}=\tilde{M}-\tilde{L}$ and note that ${\rm rank}(\tilde{L})\leq {\rm
rank}(\bar{u}\otimes 1_{n})=1$. Thus, by (i) of Lemma
\ref{lemma:ESD_comparison}, we have
$\|F^{\tilde{B}\tilde{B}'}-F^{\tilde{M}\tilde{M}'}\|\leq (n-q)^{-1}$; hence,
the ESD of $\tilde{B}\tilde{B}'$ converges a.s. in distribution to the M-P law,
and $\lambda_{\min}(\tilde{B}\tilde{B}')$ and $\lambda_{\max}(
\tilde{B}\tilde{B}')$ converge a.s. to $b_{-}(\tau)$ and $b_{+}(\tau)$,
respectively.

Next, write $\tilde{A}=(\tilde{M}-\tilde{L})D_{s}^{-1}$. By (ii) of Lemma
\ref{lemma:ESD_comparison}, we have
$L^{4}(F^{\tilde{A}\tilde{A}'},F^{\tilde{B}\tilde{B}'}) \leq
2(n-q)^{-2}(\|\tilde{A}\|_{2}^{2}+\|\tilde{B}\|_{2}^{2})\|\tilde{A}
-\tilde{B}\|_{2}^{2}$. Note that $\|\tilde{B}\|_{2}^{2}={\rm tr}(\tilde{B}
\tilde{B}')={\rm tr}(\tilde{M}\tilde{M}')-2{\rm tr}(\tilde{L}\tilde{M}') +{\rm
tr}(\tilde{L}\tilde{L}')$. By Lemma \ref{lemma:MP_law}, we have ${\rm
tr}(\tilde{M} \tilde{M}')={\rm tr}(p^{-1}MM')=(n-q)O_{\rm a.s.}(1)$, where
$O_{\rm a.s.}(1)$ denotes a term that is bounded almost surely. We have
\begin{eqnarray*}
|{\rm tr}(\tilde{L}\tilde{M}')|&=&\frac{1}{n}|{\rm tr}\{(1_{n}'p^{-1}UU'
\otimes 1_{n})(AA')\}|\\
&=&\frac{1}{n}|{\rm tr}(1_{n}'p^{-1}UU'AA'\otimes 1_{n})|\\
&=&\frac{1}{n}|1_{n}'p^{-1}UU'AA'1_{n}|\\
&\leq&\frac{1}{n}\sqrt{1_{n}'(p^{-1}UU')^{2}1_{n}}\sqrt{1_{n}'(AA')^{2}
1_{n}}\\
&\leq&\lambda_{\max}(p^{-1}UU')\lambda_{\max}(AA')\\
&=&\lambda_{\max}(p^{-1}UU'),
\end{eqnarray*}
which is $O_{\rm a.s.}(1)$ by Lemma \ref{lemma:extreme_eigenvalues}, and ${\rm
tr}(\tilde{L}\tilde{L}')\leq\lambda_{\max}(p^{-1}UU')=O_{\rm a.s.}(1)$. It
follows that $\|\tilde{B}\|_{2}^{2}=(n-q)O_{\rm a.s.}(1)$. Also, we have
$\|\tilde{A}\|_{2}^{2}={\rm tr}(\tilde{B}D_{s}^{-2}\tilde
{B}')\leq\lambda_{\max}(D_{s}^{-2})\|\tilde{B}\|_{2}^{2}=\|\tilde{B}
\|_{2}^{2}/\min_{1\leq j\leq p}s_{j}^{2}$. By Lemma
\ref{lemma:columnwise_mean_dev}, we have $\max_{1 \leq j\leq
p}|s_{j}^{2}-1|\stackrel{\rm a.s.}{\longrightarrow}0$, hence, we have
$(\min_{1\leq j\leq p}s_{j}^{2})^{-1}=O_{\rm a.s.}(1)$. It follows that
$\|\tilde{A}\|_{2}^{2}=(n-q)O_{\rm a.s.}(1)$. Finally, we have $\|\tilde
{A}-\tilde{B}\|_{2}^{2}={\rm tr}\{\tilde{B}(I_{p}-D_{s}^{-1})^{2}\tilde
{B}'\}\leq\lambda_{\max}\{(I_{p}-D_{s}^{-1})^{2}\}\|\tilde{B}\|_{2}^{2}$, and
\begin{eqnarray}
\lambda_{\max}\{(I_{p}-D_{s}^{-1})^{2}\}&\leq&\frac{(\max_{1\leq
j\leq p}|s_{j}^{2}-1|)^{2}}{(\min_{1\leq j\leq p}s_{j}^{2}+\min_{1\leq
j\leq p}s_{j})^{2}}\nonumber\\
&=&o_{\rm a.s.}(1).
\end{eqnarray}
It follows that $\|\tilde{A}-\tilde{B}\|_{2}^{2}=(n-q)o_{\rm a.s.}(1)$.
Thus, we have $L^{4}(F^{\tilde{A}\tilde{A}'},F^{\tilde{B}\tilde{B}'})
=o_{\rm a.s.}(1)$, hence the ESD of $\tilde{A}\tilde{A}'$ converges a.s.
in distribution to the M-P law.

Note that $\tilde{A}\tilde{A}'=\tilde{B}\tilde{B}'+\Delta$ with $\Delta
=\tilde{B}(D_{s}^{-2}-I_{p})\tilde{B}$, hence $\lambda_{\max}(\tilde
{A}\tilde{A}')\geq\lambda_{\max}(\tilde{B}\tilde{B}')-\|\Delta\|$ and
$\lambda_{\max}(\tilde{A}\tilde{A}')\leq\lambda_{\max}(\tilde{B}\tilde
{B}')+\|\Delta\|$ (e.g., Jiang 2010, p. 167; also using the fact that
$\lambda_{\max}(M)\leq\|M\|$ and $\lambda_{\min}(M)\geq-\|M\|$ for any
symmetric matrix $M$). Similarly, we have $\lambda_{\min}(\tilde{A}\tilde
{A}')\geq\lambda_{\min}(\tilde{B}\tilde{B}')-\|\Delta\|$ and
$\lambda_{\min}(\tilde{A}\tilde{A}')\leq\lambda_{\min}(\tilde{B}\tilde
{B}')+\|\Delta\|$. It remains to show that $\|\Delta\|\stackrel{\rm
a.s.}{\longrightarrow}0$, but this follows from
$$
\|\Delta\|\leq\|\tilde{B}\|^{2}\|D_{s}^{-2}-I_{p}\|\\
\leq\frac{\max_{1\leq j\leq p}|s_{j}^{2}-1|}{\min_{1\leq j\leq p}
s_{j}^{2}}\lambda_{\max}(\tilde{B}\tilde{B}')\stackrel{\rm
a.s.}{\longrightarrow}0.
$$

\subsection{Notation}\label{subsec:notation}
\hspace{4mm}

Some notation will be used throughout the next two subsections. Most of
these have been introduced before; we summarize below for convenience.
Recall that $A$ is an $n\times(n-q)$ matrix with $A'X=0$ and $A'A=I_{n-q}$.
We write $Z=[Z_{(1)}\;Z_{(2)}]$, where $Z_{(1)}$ is $n\times m$ and
$Z_{(2)}$ is $n\times(p-m)$, $\tilde{Z}=p^{-1/2}Z$, and $\tilde{Z}_{(j)}
=p^{-1/2}Z_{(j)}$, $j=1,2$. Also, we have $\tilde{y}=y-X\beta=
\tilde{Z}_{(1)}\alpha_{(1)}+\epsilon$ so that $\tilde{y}|Z\sim N(0,
\sigma_{\epsilon 0}^{2}V_{1,0})$, where $V_{1,0}=I_{n}+\gamma_{0}
\tilde{Z}_{(1)}\tilde{Z}_{(1)}'$; similarly, $\Sigma_{1,0}=A'V_{1,0}A=
I_{n-q}+\gamma_{0}A'\tilde{Z}_{(1)}\tilde{Z}_{(1)}'A$. Moreover, let
$\zeta=A'Z$; $\bar{U}=p^{-1}\zeta\zeta'=\tilde{\zeta}\tilde{\zeta}'$ with
$\tilde{\zeta}=p^{-1/2}\zeta$; $V_{\gamma}=I_{n}+\gamma\tilde{Z}\tilde{Z}'
=I_{n}+(\gamma/p)ZZ'$; $P_{\gamma}=A\Sigma_{\gamma}^{-1}A'=V_{\gamma}^{-1}
-V_{\gamma}^{-1}X(X'V_{\gamma}^{-1}X)^{-1}X'V_{\gamma}^{-1}$ with
$\Sigma=\Sigma_{\gamma}=I_{n-q}+\gamma\bar{U}$ (e.g., Jiang 2007, p. 13);
$G=G_{\gamma}=-(\partial/\partial\gamma)\Sigma_{\gamma}^{-1}
=\Sigma_{\gamma}^{-1}\bar{U}\Sigma_{\gamma}^{-1}$. Define $b_{1}(\gamma)
={\rm tr}(\Sigma_{\gamma}^{-1}\bar{U}\Sigma_{\gamma}^{-1}\Sigma_{1,0})$,
$b_{2}(\gamma)={\rm tr}(\Sigma_{\gamma}^{-2}\Sigma_{1,0})$, $c_{1}(\gamma)
={\rm tr}(\Sigma_{\gamma}^{-1}\bar{U})$, $c_{2}(\gamma)={\rm
tr}(\Sigma_{\gamma}^{-1})$, $s(\gamma)=y'P_{\gamma}^{2}y/{\rm
tr}(P_{\gamma})=\tilde{y}P_{\gamma}^{2}\tilde{y}/{\rm tr}(P_{\gamma})$, and
$\Delta(\gamma)=y'B_{\gamma}y=\tilde{y}'B_{\gamma}\tilde{y}$ with
$$
B=B_{\gamma}=\frac{P_{\gamma}\tilde{Z}\tilde{Z}'P_{\gamma}}{{\rm
tr}(P_{\gamma}\tilde{Z}\tilde{Z}')}-\frac{P_{\gamma}^{2}}{{\rm
tr}(P_{\gamma})}.
$$
Finally, we introduce the function
\begin{eqnarray}\label{eq:h_k_l}
h_{k,l}(\gamma)&=&\int\frac{x^{l}}{(1+\gamma x)^{k}}f_{\tau}(x)dx,
\end{eqnarray}
where $f_{\tau}$ denotes the pdf of the M-P law with the parameter $\tau\in
(0,1]$. Some special cases are, with the notation,
$$f_{1}(\gamma)=h_{2,1}(\gamma),\;f_{2}(\gamma)=h_{2,0}(\gamma),\;
g_{1}(\gamma)=h_{1,1}(\gamma),\;g_{2}(\gamma)=h_{1,0}(\gamma),$$
We shall also write $\gamma_{*}=\omega\gamma_{0}$.

\subsection{Proof of Theorem
\ref{thm:consistency_REML}}\label{subsec:proof_thm_consistency}
\hspace{4mm}

Our approach is to first consider a simplified version of Theorem
\ref{thm:consistency_REML}, in which the entries of $Z$ are i.i.d. $N(0,1)$,
and then extend the proof by explaining how to relax the restriction.

{\it Part (i).} First consider the asymptotic hehavior of $\hat{\gamma}$. For
any fixed $\gamma>0$, write $\Delta=\Delta(\gamma)$ and $B=B_{\gamma}$ for
notational simplicity. Note that $\zeta$ is $(n-q)\times p$, whose entries are
independent $N(0,1)$. Straight calculation, and Corollary
\ref{cor:MP_corollary}, show that ${\rm tr}(P_{\gamma})={\rm
tr}(\Sigma^{-1})=O_{\rm P}(n)$, and ${\rm
tr}(P_{\gamma}\tilde{Z}\tilde{Z}')={\rm tr}(\Sigma^{-1}\bar{U}) =O_{\rm P}(n)$.

Next, write $\Delta={\rm E}(\Delta|Z)+\Delta-{\rm E}(\Delta|Z)=\Delta_{1}
+\Delta_{2}$. By the normal theory (e.g., Jiang 2007, p. 238), it can be shown
that ${\rm var}(\Delta|Z)=2\sigma_{\epsilon 0}^{2}{\rm tr}(DV_{1,0}DV_{1,0})$,
where $D=A(C_{1}-C_{2})A'$ with $C_{1}= \Sigma^{-1}\bar{U}\Sigma^{-1}/c_{1}$,
$c_{1}={\rm tr}(\Sigma^{-1}\bar{U})$, $C_{2}=\Sigma^{-2}/c_{2}$, and
$c_{2}={\rm tr}(\Sigma^{-1})$. By Corollary \ref{cor:MP_corollary}, we have
$c_{j}=O_{\rm P}(n), j=1,2$. On the other hand, we have ${\rm
tr}(DV_{1,0}DV_{1,0})={\rm tr}[\{(C_{1}-C_{2})\Sigma_{1,0}\}^{2} ]={\rm
tr}\{(C_{1}\Sigma_{1,0})^{2}\}-2{\rm tr}(C_{1}\Sigma_{1,0}C_{2}
\Sigma_{1,0})+{\rm tr}\{(C_{2}\Sigma_{1,0})^{2}\}$;
$$
{\rm tr}\{(C_{1}\Sigma_{1,0})^{2}\}\leq c_{1}^{-2}{\rm
tr}\{(\Sigma^{-1}\bar{U}\Sigma^{-1}\Sigma_{0})^{2}\}=O_{\rm P}(n^{-1}),
$$
by Corollary \ref{cor:MP_corollary}, where $\Sigma_{0}$ is $\Sigma$ with
$\gamma$ replaced by $\gamma_{0}$; $${\rm tr}\{(C_{2}\Sigma_{1,0})^{2}\}\leq
c_{2}^{-2}{\rm tr}\{(\Sigma^{-1}\Sigma_{0}\Sigma^{-1})^{2}\}=O_{\rm
P}(n^{-1}),\;{\rm and}$$
$$
{\rm tr}(C_{1}\Sigma_{1,0}C_{2}\Sigma_{1,0})=
(c_{1}c_{2})^{-1}{\rm tr}(\Sigma^{-1}\bar{U}\Sigma^{-1}\Sigma_{1,0}
\Sigma^{-2}\Sigma_{1,0})=O_{\rm P}(n^{-1}).
$$
It follows that ${\rm var}(\Delta|Z)=O_{\rm P}(n^{-1})$, hence, for any
$\delta>0$, we have ${\rm P}\{|\Delta-{\rm E}(\Delta|Z)|>\delta|Z\}\leq
\delta^{-2}{\rm var}(\Delta|Z)\stackrel{\rm P}{\longrightarrow}0$, as
$n\rightarrow\infty$. Thus, by the dominated convergence theorem, we have
${\rm P}\{|\Delta-{\rm E}(\Delta|Z)|>\delta\}\rightarrow 0$, $\forall
\delta>0$, implying $\Delta_{2}=o_{\rm P}(1)$.

Next, we have $\Delta_{1}={\rm E}(\Delta|Z)=\sigma_{\epsilon 0}^{2}
(b_{1}/c_{1}-b_{2}/c_{2})$, $b_{1}={\rm tr}(\Sigma^{-1}\bar{U}
\Sigma^{-1}\Sigma_{1,0})$, $b_{2}={\rm tr}(\Sigma^{-2}\Sigma_{1,0})$, and
$c_{1}, c_{2}$ are defined earlier. By Lemma \ref{lemma:MP_law}, we have
\begin{eqnarray}\label{eq:c1}
\frac{c_{1}}{n-q}=\frac{1}{n-q}\sum_{k=1}^{n-q}\frac{\lambda_{k}}{1+
\gamma\lambda_{k}}\stackrel{\rm a.s.}{\longrightarrow}
\int_{b_{-}(\tau)}^{b_{+}(\tau)}\frac{xf_{\tau}(x)}{1+\gamma x}dx,
\end{eqnarray}
where $\lambda_{k}, 1\leq k\leq n-q$ are the eigenvalues of $\bar{U}$.
Similarly, we have
\begin{eqnarray}\label{eq:c2}
\frac{c_{2}}{n-q}=\frac{1}{n-q}\sum_{k=1}^{n-q}\frac{1}{1+\gamma
\lambda_{k}}\stackrel{\rm a.s.}{\longrightarrow}
\int_{b_{-}(\tau)}^{b_{+}(\tau)}\frac{f_{\tau}(x)}{1+\gamma x}dx.
\end{eqnarray}
Also, we have $b_{1}={\rm tr}(\Sigma^{-1}\bar{U}\Sigma^{-1})+
\gamma_{0}{\rm tr}\{\Sigma^{-1}\bar{U}\Sigma^{-1}\bar{U}_{(1)}\}$, and
\begin{eqnarray}\label{eq:tr_U_bar}
\frac{{\rm tr}(\Sigma^{-1}\bar{U}\Sigma^{-1})}{n-q}=\frac{1}{n-q}
\sum_{k=1}^{n-q}\frac{\lambda_{k}}{(1+\gamma\lambda_{k})^{2}}
\stackrel{\rm a.s.}{\longrightarrow}\int_{b_{-}(\tau)}^{b_{+}(\tau)}
\frac{xf_{\tau}(x)}{(1+\gamma x)^{2}}dx.
\end{eqnarray}
On the other hand, note that ${\rm tr}\{\Sigma^{-1}\bar{U}\Sigma^{-1}
\bar{U}_{(1)}\}=p^{-1}\sum_{k=1}^{m}\zeta_{k}'G\zeta_{k}$, where $\zeta_{k}$ is
the $k$-th column of $\zeta$, and $G=\Sigma^{-1}\bar{U} \Sigma^{-1}$. Write
$\Sigma=\Sigma_{-k}+(\gamma/p)\zeta_{k}\zeta_{k}'$, where
$\Sigma_{-k}=I_{n-q}+(\gamma/p)\sum_{l\neq k}\zeta_{l}\zeta_{l}'$. Using a
matrix identity (e.g., Sen and Srivastava 1990, p. 275), we have
$\Sigma^{-1}=\Sigma_{-k}^{-1}-(\gamma/p)\{1+(\gamma/p)
u_{k}\}^{-1}\Sigma_{-k}^{-1}\zeta_{k}\zeta_{k}'\Sigma_{-k}^{-1}$, where
$u_{k}=\zeta_{k}'\Sigma_{-k}^{-1}\zeta_{k}$. Thus, after some tedious
derivation, we have the expression
\begin{eqnarray}\label{eq:G_zeta_quad}
\zeta_{k}'G\zeta_{k}&=&\frac{u_{k}^{2}+v_{k}}{p\{1+(\gamma/p)u_{k}
\}^{2}},
\end{eqnarray}
where $v_{k}=\zeta_{k}'\Sigma_{-k}^{-1}U_{-k}\Sigma_{-k}^{-1}\zeta_{k}$ and
$U_{-k}=\sum_{l\neq k}\zeta_{l}\zeta_{l}'$. Note that $\zeta_{k}$ is
independent with $\Sigma_{-k}$. Thus, by Proposition \ref{prop:Hanson_Wright},
we have, for any $1\leq k\leq m$ and $t>0$,
\begin{eqnarray}\label{eq:u_k_concentration}
&&{\rm P}\{|u_{k}-{\rm tr}(\Sigma_{-k}^{-1})|>t|\Sigma_{-k}\}\nonumber\\
&\leq&2\exp\left\{-c\min\left(\frac{t^{2}}{K^{4}\|\Sigma_{-k}^{-1}
\|_{2}^{2}},\frac{t}{K^{2}\|\Sigma_{-k}^{-1}\|}\right)\right\},
\end{eqnarray}
where $c$ and $K$ are some positive constants. If we let
$$
t=t_{m,k}=K^{2}\max\left(\sqrt{\frac{2\log(m)}{c}}\|\Sigma_{-k}^{-1}
\|_{2},\frac{2\log(m)}{c}\|\Sigma_{-k}^{-1}\|\right),
$$
then, it is seen that the $\min$ in (\ref{eq:u_k_concentration}) is $\geq
2\log(m)/c$. It follows that ${\rm P}\{t_{m,k}^{-1}|u_{k}-{\rm
tr}(\Sigma_{-k}^{-1})|>1 |\Sigma_{-k}\}\leq 2/m^{2}, 1\leq k\leq m$, hence
\begin{eqnarray}\label{eq:u_k_prob_bound}
{\rm P}\left\{\max_{1\leq k\leq m}t_{m,k}^{-1}|u_{k}-{\rm
tr}(\Sigma_{-k}^{-1})|>1\right\}&\leq&\frac{2}{m}.
\end{eqnarray}

On the other hand, we have $\|\Sigma_{-k}^{-1}\|\leq 1$, and
$\|\Sigma_{-k}^{-1}\|_{2}\leq\sqrt{{\rm tr}(\Sigma^{-2})+8}=O_{\rm
P}(\sqrt{n})$, by Corollary \ref{cor:MP_corollary}. It follows by
(\ref{eq:m_n_p_condition}) that
\begin{eqnarray}\label{eq:u_k_order}
\max_{1\leq k\leq m}|u_{k}-{\rm tr}(\Sigma_{-k}^{-1})|&=&O_{\rm
P}(\sqrt{n\log n}).
\end{eqnarray}
Similarly, write $V_{k}=\Sigma_{-k}^{-1}U_{-k}\Sigma_{-k}^{-1}$. By a
similar argument, it can be shown that
\begin{eqnarray}\label{eq:v_k_order}
\max_{1\leq k\leq m}|v_{k}-{\rm tr}(V_{k})|&=&O_{\rm P}(n\sqrt{n\log n}).
\end{eqnarray}
Also, by an earlier expansion, it can be shown that
\begin{eqnarray}\label{eq:Sigma_k_perturbation}
|{\rm tr}(\Sigma_{-k}^{-1})-{\rm tr}(\Sigma^{-1})|&=&\frac{(\gamma/p)
\zeta_{k}'\Sigma_{-k}^{-2}\zeta_{k}}{1+(\gamma/p)u_{k}}\;\;\leq\;\;1.
\end{eqnarray}
It follows, by (\ref{eq:u_k_order}) and (\ref{eq:Sigma_k_perturbation}), that
\begin{eqnarray}\label{eq:u_k_Sigma_order}
\max_{1\leq k\leq m}|u_{k}-{\rm tr}(\Sigma^{-1})|&=&O_{\rm P}(\sqrt{n\log
n}).
\end{eqnarray}
Furthermore, by the same expansion, and (\ref{eq:Sigma_k_perturbation}), it can
be shown that
\begin{eqnarray}\label{eq:tr_V_k_order}
|{\rm tr}(V_{k})-{\rm tr}(\Sigma^{-1}U\Sigma^{-1})|&\leq&8p
\lambda_{\max}(\bar{U})+(1+2\sqrt{2})u_{k}\;\;\leq\;\;O_{\rm P}(n),
\end{eqnarray}
where the $O_{\rm P}$ does not depend on $k$. It follows, by
(\ref{eq:v_k_order}) and (\ref{eq:tr_V_k_order}), that
\begin{eqnarray}\label{eq:v_k_Sigma_order}
\max_{1\leq k\leq m}|v_{k}-p{\rm tr}(\Sigma^{-1}\bar{U}\Sigma^{-1})|
&=&O_{\rm P}(n\sqrt{n\log n}).
\end{eqnarray}
By (\ref{eq:G_zeta_quad}), (\ref{eq:u_k_Sigma_order}), and
(\ref{eq:v_k_Sigma_order}), it can be shown that $a_{1}-O_{\rm P}(\sqrt {\log
n/n})<\zeta_{k}'G\zeta_{k}/(n-q)<a_{1}+O_{\rm P}(\sqrt{\log n/n})$, where the
$O_{\rm P}$s do not depend on $k$, and
\begin{eqnarray*}
a_{1}&=&\left\{1+\gamma\left(\frac{n-q}{p}
\right)\frac{{\rm tr}(\Sigma^{-1})}{n-q}\right\}^{-2}\\
&&\times\left[\left(\frac
{n-q}{p}\right)\left\{\frac{{\rm tr}(\Sigma^{-1})}{n-q}\right\}^{2}+
\frac{{\rm tr}(\Sigma^{-1}\bar{U}\Sigma^{-1})}{n-q}\right],
\end{eqnarray*}
It then follows, by Lemma \ref{lemma:MP_law}, that ${\rm
tr}\{\Sigma^{-1}\bar{U} \Sigma^{-1}\bar{U}_{(1)}\}/(n-q)\stackrel{\rm
P}{\longrightarrow}\omega d_{1}$, where
\begin{eqnarray*}
d_{1}&=&\left\{1+\gamma\tau\int_{b_{-}(\tau)}^{b_{+}(\tau)}\frac
{f_{\tau}(x)}{1+\gamma x}dx\right\}^{-2}\\
&&\times\left[\tau\left\{\int_{b_{-}(\tau)}^{b_{+}(\tau)}
\frac{f_{\tau}(x)}{1+\gamma x}dx\right\}^{2}+\int_{b_{-}(\tau)}^{b_{+}(\tau)}
\frac{xf_{\tau}(x)}{(1+\gamma x)^{2}}dx\right].
\end{eqnarray*}
Therefore, we have $b_{1}/c_{1}\stackrel{\rm P}{\longrightarrow}(f_{1}+
\gamma_{0}\omega d_{1})/g_{1}$.

By a similar argument, we have $b_{2}/c_{2}\stackrel{\rm
P}{\longrightarrow}(f_{2}+\gamma_{0}\omega d_{2})/g_{2}$, where
\begin{eqnarray*}
d_{2}=\left\{1+\gamma\tau\int_{b_{-}(\tau)}^{b_{+}(\tau)}\frac
{f_{\tau}(x)}{1+\gamma x}dx\right\}^{-2}\int_{b_{-}(\tau)}^{b_{+}(\tau)}
\frac{f_{\tau}(x)}{(1+\gamma x)^{2}}dx.
\end{eqnarray*}

We have proved that $\Delta_{1}$ converges in probability to a constant
limit. The next thing we do is to determine the limit, in a different
way. This is because the expression of the limit given above involving
the $d$'s is a bit complicated, from which it is not easy to make a
conclusion. To this end, it is easy to show that $0\leq b_{j}/c_{j}\leq
(\gamma_{0}/\gamma)\vee 1, j=1,2$. Thus, by the dominated convergence
theorem, ${\rm E}(b_{j}/c_{j})$ converges to the same limit as
$b_{j}/c_{j}$, $j=1,2$. On the other hand, it can be shown that
\begin{eqnarray}
{\rm E}\left(\frac{b_{1}}{c_{1}}\right)&=&{\rm E}\left\{\frac{{\rm
tr}(G)}{c_{1}}\right\}+\gamma_{0}\left(\frac{m}{p}\right){\rm
E}\left\{\frac{{\rm tr}(G\bar{U})}{c_{1}}\right\}, \label{eq:b1_c1_ratio} \\
{\rm E}\left(\frac{b_{2}}{c_{2}}\right)&=&{\rm E}\left\{\frac{{\rm
tr}(\Sigma^{-2})}{c_{2}}\right\}+\gamma_{0}\left(\frac{m}{p}\right){\rm
E}\left\{\frac{{\rm tr}(\Sigma^{-2}\bar{U})}{c_{2}}\right\}.
\label{eq:b2_c2_ratio}
\end{eqnarray}
Furthermore, it is easy to show that $0\leq{\rm tr}(G)/c_{1}\leq 1$, $0\leq{\rm
tr}(G\bar{U})/c_{1}\leq\gamma^{-1}$, $0\leq{\rm tr}(\Sigma^{-2}) /c_{2}\leq 1$,
and $0\leq{\rm tr}(\Sigma^{-2}\bar{U})/c_{2}\leq\gamma^{-1}$. Thus, by Lemma
\ref{lemma:MP_law} and, again, the dominated convergence theorem, the right
sides of (\ref{eq:b1_c1_ratio}) and (\ref{eq:b2_c2_ratio}) converge to the
limit $l_{1}, l_{2}$, respectively, where $l_{j}=u_{j}+\gamma_{0}\omega w_{j},
u_{j}=f_{j}/g_{j}, j=1,2$, $w_{1}=
\int_{b_{-}(\tau)}^{b_{+}(\tau)}\{x^{2}f_{\tau}(x)/(1+\gamma x)^{2}\}dx
/g_{1}$, and $w_{2}=f_{1}/g_{2}$. Thus, with a little bit of algebra, it
follows that the limit of $\Delta_{1}$ is $\sigma_{\epsilon 0}^{2}\{
(\gamma_{*}/\gamma)-1\}(u_{2}-u_{1})$, and $u_{2}-u_{1}>0$ by a well-known
inequality (e.g., Jiang 2010, pp. 147-148).

Finally, recall that $\Delta=\Delta(\gamma)$. Thus, in conclusion, we
have shown that $\Delta(\gamma)$ converges in probability to a constant
limit, which is $>0$, $=0$, or $<0$ depending on whether $\gamma$ is
$<\gamma_{*}$, $=\gamma_{*}$, or $>\gamma_{*}$. This proves (i).

\vskip.1in {\it Part (ii).} Write $\xi=A'\tilde{y}$. We have
\begin{eqnarray*}
s'(\gamma)&=&\frac{{\rm tr}(G)\xi'\Sigma^{-2}\xi}{\{{\rm
tr}(\Sigma^{-1})\}^{2}}-2\frac{\xi'G\Sigma^{-1}\xi}{{\rm tr}(\Sigma^{-1})},
\end{eqnarray*}
It is easy to show that ${\rm E}(\xi'\xi)\leq\sigma_{\epsilon 0}^{2}(1
+\gamma_{0})(n-q)$. Thus, we have $0\leq\xi'\Sigma^{-2}\xi\leq\xi'\xi =O_{\rm
P}(n-q)$, $0\leq\xi'G\Sigma^{-1}\xi\leq\lambda_{\max}(\bar{U}) \xi'\xi=O_{\rm
P}(n-q)$, by Lemma \ref{lemma:extreme_eigenvalues}, and ${\rm
tr}(G)\leq\lambda_{\max}( \bar{U})(n-q)O_{\rm P}(n-q)$. Furthermore, for any
$0<\gamma\leq 2\gamma_{0}$, we have $(n-q)^{-1}{\rm
tr}(\Sigma^{-1})\geq(n-q)^{-1}{\rm
tr}\{(I_{n-q}+2\gamma_{0}\bar{U})^{-1}\}\stackrel{\rm a.s.}{\longrightarrow}
\int_{b_{-}(\tau)}^{b_{+}(\tau)}\{1+2\gamma_{0}x\}^{-1}f_{\tau}(x)dx>0$, by
Lemma \ref{lemma:MP_law}. Note that the $O_{\rm P}$'s here do not depend on
$\gamma$. It follows that $\sup_{0<\gamma\leq 2\gamma_{0}}|s'(\gamma)|= O_{\rm
P}(1)$. Therefore, by the Taylor expansion, we have $\hat
{\sigma}_{\epsilon}^{2}=s(\gamma_{*})+s'(\tilde{\gamma})(\hat{\gamma}
-\gamma_{*})=s(\gamma_{*})+o_{\rm P}(1)$, by part (i) of Theorem
\ref{thm:consistency_REML}, where $\tilde{\gamma}$ lies between $\gamma_{*}$
and $\hat{\gamma}$.

Next, by the proof of part (i), it is easy to show that, with $\gamma=
\gamma_{*}$, we have $s(\gamma)=\sigma_{\epsilon 0}^{2}(b_{2}/c_{2})
+O_{\rm P}(n^{-1/2})$, and $b_{2}/c_{2}\stackrel{\rm
P}{\longrightarrow}l_{2}$, where $l_{2}$ is defined in the
proof of part (i) with $\gamma=\gamma_{*}$. It follows that $l_{2}=
u_{2}+\gamma w_{2}=(f_{2}+\gamma f_{1})/g_{2}=1$. This proves part (ii).

We have proved the theorem under the assumption that the entries of $Z$ are
independent $N(0,1)$. We now explain how the result can be extended under more
general conditions. The first extension is to the case where the entries of $Z$
are i.i.d. sub-Gaussian. The only place in the proof where the normality was
used was in the early going of part (i), where the normality of $Z$ implied
that the entries of $\zeta=A'Z$ are also independent $N(0,1)$. However, the way
$A$ is involved is always through $AA'=P_{X^{\perp}}=I-P_{X}$, where
$P_{X}=X(X'X)^{-1}X'$, and $P_{X}$ has rank $q$, which is fixed (see the
beginning of Section \ref{subsec:outline_theoretical}). It turns out that
$P_{X}$ is negligible in the sense that the difference, after replacing $AA'$
by $I$, the ($n\times n$) identity matrix, it does not affect the order of the
approximation in every single place throughout the proof. Furthermore, when $A$
is replaced by $I$, the entries of $\zeta$ are clearly i.i.d., and the rest of
the proof applies without any change to the case where the entries of $Z$ are
independent sub-Gaussian. This extends the result to the latter case.

The next extension is to the case of standardized design matrix. Using the
preliminary results, namely, Lemma \ref{lemma:ESD_comparison}, Lemma
\ref{lemma:columnwise_mean_dev} and Corollary \ref{cor:ESD_psi_convergence}, it
can be shown that, the difference induced by the standardization is negligible
in the same sense.

All the extensions have been verified, step-by-step, throughout the proof to
make sure that the results of Theorem \ref{thm:consistency_REML} remain valid
for the case where $Z$ is the standardized design matrix as described in
Section \ref{subsec:outline_theoretical} (also above Corollary
\ref{cor:ESD_psi_convergence}), where the entries of $U$ are i.i.d.
sub-Gaussian. The detailed verifications, which are tedious, are omitted.

\subsection{Proof of Theorem \ref{thm:rate_REML}}\label{subsec:proof_thm_rate}
\hspace{4mm}

Recall that $\hat\gamma$ solves equation (\ref{eq:gamma_hat_equation}), and
$\hat{\sigma}_{\epsilon}^{2}$ is given by the right side of
(\ref{eq:sigma_epsilon}) with $\gamma=\hat{\gamma}$. It follows that
$\Delta(\hat{\gamma})=0$ and $\hat{\sigma}_{\epsilon}^{2}=s(\hat{\gamma})$.
Theorem \ref{thm:consistency_REML} has established that
$\hat{\gamma}\stackrel{\rm P}{\longrightarrow}\gamma_{*}$. Because
$\Delta(\hat{\gamma})=0$, by the Taylor series expansion, and some algebra, we
have
\begin{eqnarray}\label{eq:hat_gamma_diff}
\hat{\gamma}-\gamma_{*}&=&-\frac{\Delta(\gamma_{*})}{\Delta'(\gamma_{*})}
+O_{\rm P}(|\Delta(\gamma_{*})|^{2}).
\end{eqnarray}
Here we also use the fact that $\Delta'(\gamma_{*})$ converges in probability
to a nonzero quantity. Indeed, from the proof of Theorem
\ref{thm:consistency_REML}, it can be checked that $\Delta'(\gamma)$ converges
in probability, for every fixed $\gamma$, to $\Delta_{\infty}'(\gamma)$, where
\begin{eqnarray*}
\Delta_{\infty}(\gamma)&=&\sigma_{\epsilon 0}^{2}\left(\frac
{\gamma_{*}}{\gamma}-1\right)\left\{\frac{f_{2}(\gamma)}{g_{2}(\gamma)}
-\frac{f_{1}(\gamma)}{g_{1}(\gamma)}\right\},
\end{eqnarray*}
and the difference within the $\{\cdots\}$ is positive. It follows that
\begin{eqnarray}\label{eq:delta_prime}
\Delta_{\infty}'(\gamma_{*})&=&-\frac{\sigma_{\epsilon 0}^{2}}{\gamma_{*}}
\left\{\frac{f_{2}(\gamma_{*})}{g_{2}(\gamma_{*})}-\frac
{f_{1}(\gamma_{*})}{g_{1}(\gamma_{*})}\right\}\;\;<\;\;0.
\end{eqnarray}

Next, a Taylor series expansion of $s(\gamma)$ yields
$\hat{\sigma}_{\epsilon}^{2}=s(\gamma_{*})+s'(\gamma_{*})(\hat{\gamma}-
\gamma_{*})+O(|\hat{\gamma}-\gamma_{*}|^{2})$, which, combined with
(\ref{eq:hat_gamma_diff}), leads to the expansion
\begin{eqnarray}\label{eq:sigma_epsilon_expansion}
\hat{\sigma}_{\epsilon}^{2}&=&s(\gamma_{*})-\frac{s'(\gamma_{*})}{\Delta'(
\gamma_{*})}\Delta(\gamma_{*})+O_{\rm P}(|\Delta(\gamma_{*})|^{2}).
\end{eqnarray}

Write $s(\gamma)=s_{1}(\gamma)+s_{2}(\gamma)$, where $s_{1}(\gamma)={\rm
E}\{s(\gamma)|Z\}$ and $s_{2}(\gamma)=s(\gamma)-s_{1}(\gamma)$. It was shown in
the proof of Theorem \ref{thm:consistency_REML} that
$s_{1}(\gamma)=\sigma_{\epsilon 0}^{2} \{b_{2}(\gamma)/c_{2}(\gamma)\}$. Also,
we have the expression $s_{2}(\gamma) =\tilde{w}'D_{\gamma}\tilde{w}-{\rm
tr}(D_{\gamma})$, where
$$D_{\gamma}=\sigma_{\epsilon 0}^{2}\frac{\Sigma_{1,0}^{1/2}
\Sigma_{\gamma}^{-1}\bar{U}\Sigma_{\gamma}^{-1}\Sigma_{1,0}^{2}}{{\rm
tr}(\Sigma_{\gamma}^{-1})}\;\;{\rm and}\;\;
\tilde{w}=\frac{\Sigma_{1,0}^{-1/2}A'\tilde{y}}{\sigma_{\epsilon 0}}.$$ Note
that $\tilde{w}|Z\sim N(0,I_{n-q})$. Also recall (from the proof of Theorem
\ref{thm:consistency_REML}, part (i)) that
$\Delta(\gamma)(\gamma)=\Delta_{1}(\gamma) +\Delta_{2}(\gamma)$ with
$\Delta_{1}(\gamma)=\sigma_{\epsilon 0}^{2}
\sum_{j=1}^{2}(-1)^{j-1}b_{j}(\gamma)/c_{j}(\gamma)$, and, similarly,
$\Delta_{2}(\gamma)=\tilde{w}'F_{\gamma}\tilde{w}- {\rm tr}(F_{\gamma})$, where
$F_{\gamma}=\sigma_{\epsilon 0}^{2}
\Sigma_{1,0}^{1/2}H_{\gamma}\Sigma_{1,0}^{1/2}$ with
$$
H_{\gamma}=\frac{\Sigma_{\gamma}^{-1}\bar{U}\Sigma_{\gamma}^{-1}\bar
{U}\Sigma_{\gamma}^{-1}}{{\rm tr}(\Sigma_{\gamma}^{-1}\bar{U})}-\frac
{\Sigma_{\gamma}^{-1}\bar{U}\Sigma_{\gamma}^{-1}}{{\rm
tr}(\Sigma_{\gamma}^{-1})}.
$$

As in the proof of Theorem \ref{thm:consistency_REML}, part (ii), write
$\xi=A'\tilde{y}$, and observe that
\begin{eqnarray*}
s'(\gamma)&=&\frac{{\rm tr}(\Sigma_{\gamma}^{-1}\bar{U}\Sigma_{\gamma}^{-1})
\xi'\Sigma_{\gamma}^{-2}\xi}{\{{\rm tr}(\Sigma_{\gamma}^{-1})\}^{2}}-\frac
{2\xi'\Sigma_{\gamma}^{-1}\bar{U}\Sigma_{\gamma}^{-2}\bar{U}
\Sigma_{\gamma}^{-1}\xi}{{\rm tr}(\Sigma_{\gamma}^{-1})}.
\end{eqnarray*}
We have ${\rm E}(\xi'\Sigma_{\gamma}^{-2}\xi|Z)=\sigma_{\epsilon 0}^{2}{\rm
tr}(\Sigma_{\gamma}^{-2}\Sigma_{1,0})$, and
$$
{\rm E}(\xi'\Sigma_{\gamma}^{-1}\bar{U}\Sigma_{\gamma}^{-2}\bar{U}
\Sigma_{\gamma}^{-1}\xi|Z)=\sigma_{\epsilon 0}^{2}{\rm
tr}(\Sigma_{\gamma}^{-1}\bar{U}\Sigma_{\gamma}^{-2}\bar{U}
\Sigma_{\gamma}^{-1}\Sigma_{1,0}).
$$
With these, using similar derivations to the proof of Theorem
\ref{thm:consistency_REML}, we conclude that
\begin{eqnarray}\label{eq:s_prime}
s'(\gamma)&\stackrel{\rm P}{\longrightarrow}&s_{\infty}'(\gamma)
\nonumber\\
&=&\sigma_{\epsilon 0}^{2}\left[\frac{h_{2,1}(\gamma)\{h_{2,0}(\gamma)+
\gamma_{*}h_{2,1}(\gamma)\}}{\{h_{1,0}(\gamma)\}^{2}}-\frac{2\{h_{4,
2}(\gamma)+\gamma_{*}h_{4,3}(\gamma)\}}{h_{1,0}(\gamma)}\right]
\end{eqnarray}
(see (\ref{eq:h_k_l}) for notation). Thus, going back to
(\ref{eq:sigma_epsilon_expansion}), we can write
\begin{eqnarray}\label{eq:sigma_epsilon_diff}
\hat{\sigma}_{\epsilon}^{2}-\sigma_{\epsilon 0}^{2}&=&
s_{2}(\gamma_{*})-\frac{s_{\infty}'(\gamma_{*})}{\Delta_{\infty}'(
\gamma_{*})}\Delta_{2}(\gamma_{*})\nonumber\\
&&+\sigma_{\epsilon 0}^{2}\left\{\frac{b_{2}(\gamma_{*})}{c_{2}(
\gamma_{*})}-1\right\}-\sigma_{\epsilon 0}^{2}\frac{s_{\infty}'(
\gamma_{*})}{\Delta_{\infty}'(\gamma_{*})}\left\{\frac{b_{1}(
\gamma_{*})}{c_{1}(\gamma_{*})}-\frac{b_{2}(\gamma_{*})}{c_{2}(
\gamma_{*})}\right\}\nonumber\\
&&-\frac{s'(\gamma_{*})-s_{\infty}'(\gamma_{*})}{\Delta_{\infty}'(
\gamma_{*})}\Delta(\gamma_{*})+\frac{\{\Delta'(\gamma_{*})
-\Delta_{\infty}'(\gamma_{*})\}s'(\gamma_{*})}{\Delta'(\gamma_{*})
\Delta_{\infty}'(\gamma_{*})}\Delta(\gamma_{*})\nonumber\\
&&+O_{\rm P}(|\Delta(\gamma_{*})|^{2}).
\end{eqnarray}
We shall argue that all of the terms on the right side of
(\ref{eq:sigma_epsilon_diff}) except those in the second line are $o_{\rm
P}(\sqrt{\log n/n})$, while the terms in the second line are $O_{\rm
P}(\sqrt{\log n/n})$. For the last two lines, it suffices to show that
\begin{eqnarray}\label{eq:delta_gamma_order}
\Delta(\gamma_{*})&=&O_{\rm P}\left(\sqrt{\frac{\log n}{n}}\right),
\end{eqnarray}
because $s'(\gamma_{*})-s_{\infty}'(\gamma_{*})=o_{\rm P}(1)$ and
$\Delta'(\gamma_{*})-\Delta_{\infty}'(\gamma_{*})=o_{\rm P}(1)$. Note that
(\ref{eq:delta_gamma_order}) also ensures $\hat{\gamma}-\gamma_{*}=O_{\rm
P}(\sqrt{\log n/n})$ by virtue of (\ref{eq:hat_gamma_diff}), the convergence of
$\Delta'(\gamma_{*})$ to $\Delta_{\infty}'(\gamma_{*})$, and
(\ref{eq:delta_prime}). In order to establish (\ref{eq:delta_gamma_order}), we
need the following lemma.

\begin{lemma}\label{lemma:tr_psi_bound}
Suppose that (\ref{eq:m_n_p_condition}) holds and let $\Psi=p^{-1}Z'Z$ and
$\tilde{\Sigma}_{\gamma}=I_{n} +\gamma\Psi$. Then, we have
\begin{eqnarray}\label{eq:tr_psi_bound}
{\rm tr}(\Sigma_{\gamma}^{-k}\bar{U}^{l})&=&{\rm
tr}(\tilde{\Sigma}_{\gamma}^{-k}\Psi^{l})+O_{\rm P}(1),\;\;\;\forall k\geq 1,
l\geq 0.
\end{eqnarray}
\end{lemma}

The proof of Lemma \ref{lemma:tr_psi_bound}, which is omitted, follows closely
the note regarding $AA'$ near the end of the proof of Theorem
\ref{thm:consistency_REML}. The advantage of this lemma is that, because the
entries of $Z$ are independent sub-Gaussian with mean $0$, unit variance, and
bounded fourth moments, the behavior of the trace on the right side of
(\ref{eq:tr_psi_bound}) is well studied. Indeed, we can use Theorem 9.10 of Bai
and Silverstein (2010) on the asymptotic behavior of linear spectral statistics
to claim that, for all $k\geq 1, l\geq 0$, we have
\begin{eqnarray}\label{eq:tr_Sigma_psi_concentration}
\left|\frac{1}{n}{\rm tr}(\tilde{\Sigma}_{\gamma}^{-k}\Psi^{l})-\int
\frac{x^{l}}{(1+\gamma x)^{k}}f_{n/p}(x)dx\right|&=&O_{\rm P}(n^{-1}).
\end{eqnarray}
Equation (\ref{eq:tr_Sigma_psi_concentration}), combined with
(\ref{eq:m_n_p_rate}), (\ref{eq:h_k_l}) and (\ref{eq:tr_psi_bound}), imply that
for all $k\geq 1, l\geq 0$, we have
\begin{eqnarray}\label{eq:tr_Sigma_psi_concentration_h_k_l}
\left|\frac{1}{p}{\rm tr}(\Sigma_{\gamma}^{-k}\bar{U}^{l})-\tau
h_{k,l}(\gamma)\right|&=&O_{\rm P}(n^{-1}).
\end{eqnarray}
Therefore, we have $(n-q)^{-1}c_{j}(\gamma_{*})-c_{j,\infty}(\gamma_{*})
=O_{\rm P}(n^{-1}), j=1,2$, where $c_{1,\infty}(\gamma_{*})=h_{1,
1}(\gamma_{*})$ and $c_{2,\infty}(\gamma_{*})=h_{1,0}(\gamma_{*})$.

On the other hand, by the proof of Theorem \ref{thm:consistency_REML}, part
(i), and (\ref{eq:tr_psi_bound}), we have
$(n-q)^{-1}b_{j}(\gamma_{*})-b_{j,\infty}(\gamma_{*})=O_{\rm P}(\sqrt{\log
n/n}), j=1,2$, where
\begin{eqnarray*}
b_{1,\infty}(\gamma_{*})&=&h_{2,1}(\gamma_{*})+\gamma_{*}\{1+
\tau\gamma_{*}h_{1,0}(\gamma_{*})\}^{-2}\{\tau h_{1,0}^{2}(\gamma_{*})+
h_{2,1}(\gamma_{*})\},\\
b_{2,\infty}(\gamma_{*})&=&h_{2,0}(\gamma_{*})+\gamma_{*}\{1+\tau
\gamma_{*}h_{1,0}(\gamma_{*})\}^{-2}h_{2,0}(\gamma_{*}).
\end{eqnarray*}

Also, by the proof of Theorem \ref{thm:consistency_REML}, part (ii), we have
$b_{2, \infty}(\gamma_{*})/c_{2,\infty}(\gamma_{*})=1$. Moreover, $\Delta_{1}(
\gamma_{*})\stackrel{\rm P}{\longrightarrow}0$, implying $b_{1,\infty}(
\gamma_{*})/c_{1,\infty}(\gamma_{*})-b_{2,\infty}(\gamma_{*})/c_{2,
\infty}(\gamma_{*})=0$. Thus, we conclude that $\Delta_{1}(\gamma_{*}) =O_{\rm
P}(\sqrt{\log n/n})$. On the other hand, it is seen from the proof of Theorem
\ref{thm:consistency_REML} that $\Delta_{2}(\gamma_{*})=O_{\rm P}(n^{-1/2})$.
Therefore, (\ref{eq:delta_gamma_order}) holds.

By similar arguments, it can be shown that the terms in the second line of the
right side of (\ref{eq:sigma_epsilon_diff}) are $O_{\rm P}(\sqrt{\log n/n})$.

Next, by the expressions of $s_{2}(\gamma)$, $\Delta_{2}(\gamma)$, we can write
the first line of the right side of (\ref{eq:sigma_epsilon_diff}) as
$Q=\tilde{w}'M\tilde{w}-{\rm tr}(M)$, where $M=D_{\gamma_{*}}
-\{s_{\infty}'(\gamma_{*})/\Delta_{\infty}'(\gamma_{*})\}F_{\gamma_{*}}$. We
have (e.g., Jiang 2007, p. 238) ${\rm E}(Q|Z)=0$ and ${\rm var}(Q|Z)= 2{\rm
tr}(M^{2})$. From the expressions of $D_{\gamma}$ and $F_{\gamma}$, the
following expressions can be derived:
\begin{eqnarray}
\frac{{\rm tr}(D_{\gamma}^{2})}{\sigma_{\epsilon 0}^{4}}&=&\frac {{\rm
tr}(\Sigma_{\gamma}^{-1}\bar{U}\Sigma_{\gamma}^{-1}\Sigma_{1,
0}\Sigma_{\gamma}^{-1}\bar{U}\Sigma_{\gamma}^{-1}\Sigma_{1,0})}{\{{\rm
tr}(\Sigma_{\gamma}^{-1})\}^{2}},\label{eq:tr_D_gamma_2}\\
\frac{{\rm tr}(F_{\gamma}^{2})}{\sigma_{\epsilon 0}^{4}}&=&\frac {{\rm
tr}(\Sigma_{\gamma}^{-1}\bar{U}\Sigma_{\gamma}^{-1}\bar{U}
\Sigma_{\gamma}^{-1}\Sigma_{1,0}\Sigma_{\gamma}^{-1}\bar{U}
\Sigma_{\gamma}^{-1}\bar{U}\Sigma_{\gamma}^{-1}\Sigma_{1,0})}{\{{\rm
tr}(\Sigma_{\gamma}^{-1}\bar{U})\}^{2}}\label{eq:tr_F_gamma_2}\\
&&+\frac{{\rm tr}(\Sigma_{\gamma}^{-1}\bar{U}\Sigma_{\gamma}^{-1}\Sigma_{1,
0}\Sigma_{\gamma}^{-1}\bar{U}\Sigma_{\gamma}^{-1}\Sigma_{1,0})}{\{{\rm
tr}(\Sigma_{\gamma}^{-1})\}^{2}}\nonumber\\
&&-2\frac{{\rm tr}(\Sigma_{\gamma}^{-1}\bar{U}\Sigma_{\gamma}^{-1}\bar{U}
\Sigma_{\gamma}^{-1}\Sigma_{1,0}\Sigma_{\gamma}^{-1}\bar{U}
\Sigma_{\gamma}^{-1}\Sigma_{1,0})}{{\rm tr}(\Sigma_{\gamma}^{-1}){\rm
tr}(\Sigma_{\gamma}^{-1}\bar{U})},\nonumber\\
\frac{{\rm tr}(D_{\gamma}F_{\gamma})}{\sigma_{\epsilon 0}^{4}}&=&\frac
{{\rm tr}(\Sigma_{\gamma}^{-1}\bar{U}\Sigma_{\gamma}^{-1}\bar{U}
\Sigma_{\gamma}^{-1}\Sigma_{1,0}\Sigma_{\gamma}^{-1}\bar{U}
\Sigma_{\gamma}^{-1}\Sigma_{1,0})}{{\rm tr}(\Sigma_{\gamma}^{-1}){\rm
tr}(\Sigma_{\gamma}^{-1}\bar{U})}\nonumber\\
&&-\frac{{\rm tr}(\Sigma_{\gamma}^{-1}\bar{U}\Sigma_{\gamma}^{-1}
\Sigma_{1,0}\Sigma_{\gamma}^{-1}\bar{U}\Sigma_{\gamma}^{-1}\Sigma_{1,
0})}{\{{\rm tr}(\Sigma_{\gamma}^{-1})\}^{2}}. \label{eq:tr_D_F_gamma}
\end{eqnarray}
Using the same arguments as in the proof of Theorem \ref{thm:consistency_REML},
it can be shown that when multiplied by $n$, the terms on the right sides of
(\ref{eq:tr_D_gamma_2}), (\ref{eq:tr_F_gamma_2}) and (\ref{eq:tr_D_F_gamma})
converge in probability to some constants (the derivation is tedious, and
therefore omitted). In particular, it follows, again by the dominated
convergence theorem, that the first line on the right side of
(\ref{eq:sigma_epsilon_diff}) is $o_{\rm P}(\sqrt{\log n/n})$. Therefore, by
combining the proved results, we have
$\hat{\sigma}_{\epsilon}^{2}-\sigma_{\epsilon 0}^{2}=O_{\rm P}(\sqrt{\log
n/n})$, and $\hat{\gamma}-\gamma_{*}=O_{\rm P}(\sqrt{\log n/n})$ as shown
earlier.

Finally, let $t_{1}$ denote the first two lines on the right side of
(\ref{eq:sigma_epsilon_diff}), and $t_{2}$ the last two lines. We have shown
that $t_{2}=o_{\rm P}(\sqrt{\log n/n})$ and $t_{1}=O_{\rm P}(\sqrt{\log n/n})$.
Furthermore, note that the second line on the right side of
(\ref{eq:sigma_epsilon_diff}) has zero contribution to the variance of $t_{1}$
conditioning on $Z$. Thus, the argument below (\ref{eq:tr_D_F_gamma}) has shown
that $n{\rm var}(t_{1}|Z)$ converges in probability to a constant. This
completes the proof.

\vspace{5mm}

{\bf Acknowledgement.}
The authors wish to thank Professor Iain Johnstone for
helpful discussion. The research of Jiming Jiang was partially
supported by the NSF grants DMS-0809127, SES-1121794, and the NIH
grant R01-GM085205A1. The research of of Cong Li was partially
supported by the NIH grant R01-GM59507. The research of Debashis
Paul was partially supported by the NSF grant DMS-1106690. The
research of Can Yang was partially supported by the NIH grants
R01-AA11330 and R01-DA030976. The research of Hongyu Zhao was
partially supported by the NIH grant R01-GM59507, the CTSA grant
UL1-RR024139, and the Department of Veterans Affairs (VA
Cooperative Studies Program).

\end{document}